\documentclass[12pt,reqno]{amsart}

\usepackage{amsfonts,amsmath,amsthm}
\usepackage{amssymb,epsfig}
\usepackage{enumerate} 

\usepackage{color} 
\definecolor{vert}{rgb}{0,0.6,0}

\usepackage[text={425pt,650pt},centering]{geometry}
\usepackage{graphicx}
\usepackage{epsfig}
\usepackage{tikz,esvect}
\usetikzlibrary{3d,calc,shapes}
\usepackage{pgfplots}
\pgfplotsset{compat=1.18}

\usepackage{caption}
\usepackage{color} 
\definecolor{vert}{rgb}{0,0.6,0}

\numberwithin{figure}{section}

\theoremstyle{plain}
\newtheorem{thm}{Theorem}[section]

\newtheorem{defn}{Definition}
\newtheorem{quest}{Question}

\newtheorem{ex}{Example}
\newtheorem{lem}[thm]{Lemma}
\newtheorem{cor}[thm]{Corollary}
\newtheorem{prop}[thm]{Proposition}
\theoremstyle{remark}
\newtheorem{rem}{\bf{Remark}}
\numberwithin{equation}{section}

\newcommand{\N}{\mathbb{N}}

\newcommand{\R}{\mathbb{R}}

\newcommand{\T}{\mathbb{T}}
\newcommand{\Z}{\mathbb{Z}}

\newcommand{\cH}{\mathcal{H}}
\newcommand{\cL}{\mathcal{L}}

\newcommand{\cP}{\mathcal{P}}

\newcommand{\AC}{{\rm AC\,}}

\newcommand{\Lip}{{\rm Lip\,}}

\newcommand{\gam}{\gamma}

\newcommand{\ep}{\varepsilon}

\newcommand{\sig}{\sigma}

\newcommand{\ol}{\overline}

\newcommand{\supp}{{\rm supp}\,}

\newcommand{\Div}{{\rm div}\,}

\begin{document}

\title[Nonconvex Hamilton--Jacobi equations]
{Representation formulas and large time behavior for solutions to some nonconvex Hamilton--Jacobi equations}

\author[H. V. TRAN]
{Hung Vinh Tran}

\thanks{
The work of HVT is partially supported by NSF grant DMS-2348305.
}

\address[H. V. Tran]
{
Department of Mathematics, 
University of Wisconsin Madison, Van Vleck Hall, 480 Lincoln Drive, Madison, Wisconsin 53706, USA}
\email{hung@math.wisc.edu}


\begin{abstract}
We give a new representation formula for solutions to nonconvex first-order Hamilton--Jacobi equations in the periodic setting and present some applications.
We then prove the large time behavior for solutions under some additional assumptions.
\end{abstract}

\maketitle


\section{Introduction}
In this paper, we are interested in studying properties of the viscosity solution to the following first-order Hamilton--Jacobi equation
\begin{equation}\label{eq:HJ}
    \begin{cases}
        u_t + H(x,Du)=0 \qquad &\text{ in } \T^n \times (0,\infty),\\
        u(x,0)=g(x) \qquad &\text{ on } \T^n.
    \end{cases}
\end{equation}
Here, $\T^n=\R^n/\Z^n$ is the $n$-dimensional flat torus, $H\in C^2(\T^n\times \R^n)$ is a given Hamiltonian, and $g\in C^2(\T^n)$ is a given initial data.

\medskip

Our goal is twofold. 
Firstly, we give a new representation formula to the solution $u$ to \eqref{eq:HJ}, which is a generalization of the one coming from the characteristic method.
We then use this formula to deduce the representation formula in \cite{Gomes-Mitake-Tran} in the convex setting.
We also obtain Mather measures in the nonconvex setting through the new representation formula and large time averages.
Secondly, we study the large time behavior of $u$, that is, the limit of $u(\cdot,t)$ as $t\to \infty$ in some nonconvex settings.

\subsection{Assumptions}
We list here the main assumptions that we will use.
\begin{itemize}
    \item[(A1)] We assume
    \[
    \lim_{|p|\to \infty} \min_{x\in \T^n} \left(\frac{1}{2}H(x,p)^2 + D_xH(x,p)\cdot p \right)=+\infty.
    \]
    \item[(A2)] There exists $\theta>0$ such that:  For $(x,p)\in \T^n\times \R^n$,
    \[
    D_pH(x,p)\cdot p  \geq (\theta+1) H(x,p) .
    \]
\end{itemize}

To study the properties of $u$, we consider the vanishing viscosity process: For $\ep \in (0,1)$, let $u^\ep$ solve
\begin{equation}\label{eq:HJ-ep}
    \begin{cases}
        u^\ep_t + H(x,Du^\ep)=\ep \Delta u^\ep \qquad &\text{ in } \T^n \times (0,\infty),\\
        u^\ep(x,0)=g(x) \qquad &\text{ on } \T^n.
    \end{cases}
\end{equation}

Under assumption (A1), $u^\ep$ and $u$ are globally Lipschitz thanks to the classical Bernstein method (see \cite[Chapter 1]{tranbook}).
More precisely, there exists $C>0$ depending only on $\|g\|_{C^2(\T^n)}$, $H$, and $n$ such that, for $\ep\in (0,1)$,
\begin{equation}\label{eq:Lip-bound}
\begin{cases}
    \|u^\ep_t\|_{L^\infty(\T^n\times [0,\infty))}+    \|Du^\ep\|_{L^\infty(\T^n\times [0,\infty))} \leq C,\\
     \|u_t\|_{L^\infty(\T^n\times [0,\infty))}
+    \|Du\|_{L^\infty(\T^n\times [0,\infty))}
\leq C.
\end{cases}
\end{equation}
We also have that $u^\ep$ converges to $u$ locally uniformly on $\T^n\times [0,\infty)$ as $\ep \to0$.
Assumption (A1) is always in force in this paper.
And we will only impose assumption (A2) when we study the large time behavior of $u$.
Further discussions on (A2) and its related forms are given in Section \ref{subsec:LT-A2} and Remark \ref{rem:weaker A2}.

\subsection{Main results}
The linearized operator of \eqref{eq:HJ-ep} around the solution $u^\ep$ is
\[
\phi \mapsto \cL^\ep[\phi] = \phi^\ep_t +D_pH(x,Du^\ep)\cdot D\phi - \ep \Delta \phi.
\]
Fix $T>0$.
The corresponding adjoint equation to this linearized operator is
\begin{equation}\label{eq:sigma-ep}
    \begin{cases}
        -\sigma^\ep_t -\Div(D_p H(x,Du^\ep)\sigma^\ep)=\ep \Delta \sigma^\ep \qquad &\text{ in } \T^n \times (0,T),\\
        \sigma^\ep(x,T)=\nu(x) \qquad &\text{ on } \T^n.
    \end{cases}
\end{equation}
Here, $\nu$ is a Radon probability measure.
When necessary, we write $\sigma^\ep=\sigma^{\ep,\nu}$ to specify the clear dependence.
We typically choose $\nu=\delta_z$, which is the Dirac delta measure at a given point $z\in \T^n$.
In this case, we write $\sigma^{\ep,z}=\sigma^{\ep,\delta_z}$ for simplicity.
We have $\sigma^{\ep,z}\geq 0$ on $\T^n \times [0,T)$, and
\[
\int_{\T^n}\sigma^{\ep,z}(x,t)\,dx=1 \qquad \text{ for all } t\in [0,T).
\]
As $\|Du^\ep\|_{L^\infty(\T^n\times [0,\infty))} \leq C$, for each $t\in [0,T]$, let $\nu^{\ep,z}(\cdot,\cdot,t)$ be the probability measure on $\T^n \times \ol B(0,C)$ such that, for all $\psi \in C(\T^n\times \R^n)$,
\[
\int_{\T^n\times \R^n}\psi(x,p)\,d\nu^{\ep,z}(x,p,t) = \int_{\T^n}\psi(x,Du^\ep(x,t))\sigma^{\ep,z}(x,t)\,dx.
\]
Denote by $\mu^{\ep,z}$ the measure on $\T^n\times \ol B(0,C)\times [0,T]$ such that 
\[
d\mu^{\ep,z}(x,p,t)=d\nu^{\ep,z}(x,p,t)dt.
\]
Then, $\mu^{\ep,z}$ is a nonnegative Radon measure on $\T^n\times \ol B(0,C)\times [0,T]$, and $\mu^{\ep,z}(\T^n\times \ol B(0,C)\times [0,T])=T$.
Note that $\mu^{\ep,z}, \nu^{\ep,z}, \sig^{\ep,z}$ are dependent on $T$.

\begin{thm}\label{thm:rep1}
    Assume {\rm (A1)}.
    Then, for $(z,T)\in \T^n \times (0,\infty)$, we have
    \begin{align*}
    &u^\ep(z,T) \\
    = &\int_{\T^n} g(x)\sigma^{\ep,z}(x,0)\,dx + \int_0^T\int_{\T^n\times \R^n} \left(D_pH(x,p)\cdot p - H(x,p) \right)\, d\nu^{\ep,z}(x,p,t)dt\\
     = &\int_{\T^n} g(x)\sigma^{\ep,z}(x,0)\,dx + \int_{\T^n\times \R^n\times [0,T]} \left(D_pH(x,p)\cdot p - H(x,p) \right) \,d\mu^{\ep,z}(x,p,t).
    \end{align*}
    Pick a sequence $\{\ep_k\}$ converging to $0$ such that $\sigma^{\ep_k,z}(x,0)\,dx$ converges to $d\sigma^z(x)$ and $\mu^{\ep_k,z}$ converges to $\mu^z$ weakly in the sense of measures. 
    Then,
    \begin{equation}\label{eq:formula-u-1}
          u(z,T)
     = \int_{\T^n} g(x)\,d\sigma^z(x)
     + \int_{\T^n\times \R^n\times [0,T]} \left(D_pH(x,p)\cdot p - H(x,p) \right) \,d\mu^{z}(x,p,t).  
    \end{equation}
\end{thm}

We note that \eqref{eq:formula-u-1} is new and is a generalization of the formula coming from the method of characteristics (see Remark \ref{rem:char}).
In the above, $\mu^z$ and $\sigma^z$ also depend on $T$,
and we write $\mu^z=\mu^{z,T}$, $\sigma^z=\sigma^{z,T}$ to demonstrate the clear dependence when needed.

\begin{rem}\label{rem:1}
The nonlinear adjoint method was introduced first in \cite{Evans1} to study the gradient shock structures of the Cauchy problem for nonconvex Hamiltonians. 
The static cases were studied in \cite{Tran1}.
When $H(x,p)=H(p)$, \cite{Evans2} gave a new representation formula for $u$ using generalized envelopes of affine solutions, which is different from \eqref{eq:formula-u-1}.
When $H$ depends on $x$, this is not possible as affine functions are not special solutions  to \eqref{eq:HJ}.
The approach to prove Theorem \ref{thm:rep1} is related to that in \cite{Cagnetti-Gomes-Tran}.
We refer the reader to \cite{Le-Mitake-Tran, GPV, tranbook} and the references therein for overviews of the nonlinear adjoint method.

Under (A1), there exists $C>0$ independent of $\ep\in (0,1)$ such that, for $(z,T)\in \T^n\times [0,\infty)$ and $\ep\in (0,1)$,
\[
|u^\ep(z,T)-u(z,T)| \leq C(1+T) \ep^{1/2}.
\]
See \cite{Evans1, Tran1} and the references therein.
The convergence rate $O(\ep^{1/2})$ is optimal for general Hamiltonians (see some examples in \cite{QSTY}).
We include the proof of this $O(\ep^{1/2})$ rate in Lemma \ref{lem:rate 1-2} for completeness.
If $H$ is uniformly convex in $p$, it was shown that the convergence rate is improved to $\ep|\log \ep|$ in \cite{Chaintron-Daudin, Cirant-Goffi}, which is optimal.
To be more precise, the $\ep|\log \ep|$ rate was proved for a purely quadratic Hamiltonian in one dimension in \cite{QSTY}, for general quadratic Hamiltonians with potential energies in \cite{Chaintron-Daudin}, and for uniformly convex Hamiltonians in \cite{Cirant-Goffi}.
We refer the reader to \cite{CGM23, Cirant-Goffi} for some other new quantitative convergence results via the nonlinear adjoint method.

\end{rem}

Next, we consider the situation where $H$ is convex in $p$.
Let $\mathcal{R}^+(\T^n\times \ol B(0,C)\times [0,T])$ be the set of all nonnegative Radon measures on $\T^n\times \ol B(0,C)\times [0,T]$.
Let $\cP(\T^n)$ and $\cP(\T^n \times \ol B(0,C))$ be the sets of all Radon probability measures on $\T^n$ and $\T^n \times \ol B(0,C)$, respectively.
For $\sigma_0\in \cP(\T^n)$, $\cH(\sigma_0,\delta_z;0,T)$ is the set of all $\gamma \in \mathcal{R}^+(\T^n\times \ol B(0,C)\times [0,T])$ such that
\begin{align*}
    \int_{\T^n\times \R^n\times [0,T]} (\varphi_t(x,t) + v\cdot D\varphi(x,t))\,d\gamma(x,v,t)= \varphi(z,T) - \int_{\T^n}\varphi(x,0)\,d\sigma_0(x)
\end{align*}
for all test functions $\varphi \in C^2(\T^n\times [0,T])$.
And set
\[
\cH(\delta_z;0,T) = \bigcup_{\sigma_0 \in \cP(\T^n)}\cH(\sigma_0,\delta_z;0,T).
\]
We will see below that $\cH(\delta_z;0,T) \neq \emptyset$.
For $\gamma \in \cH(\delta_z;0,T)$, we denote by $\sigma^\gamma\in\cP(\T^n)$ the unique element such that $\gamma\in \cH(\sigma^\gamma,\delta_z;0,T)$.
We are ready to state our second main result.

\begin{thm}\label{thm:rep2}
    Assume {\rm (A1)}.
    Assume further that $p\mapsto H(x,p)$ is convex for each $x\in \T^n$.
    Let $L=L(x,v)$ be the Legendre transform of $H$.
    Fix $(z,T)\in \T^n \times (0,\infty)$.
    Let $\sig^z$ and $\mu^z$ be the measures as in the statement of Theorem \ref{thm:rep1}.
    Let $\gamma^z \in \mathcal{R}^+(\T^n\times \ol B(0,C)\times [0,T])$ be such that
    \[
    \int_{\T^n\times \R^n\times [0,T]} \psi\left(x,D_pH(x,p)\right) \,d\mu^{z}(x,p,t) = \int_{\T^n\times \R^n\times [0,T]} \psi\left(x,v\right) \,d\gamma^{z}(x,v,t)
    \]
    for all $\psi\in C(\T^n\times \R^n)$.
    Then,
    \begin{equation}\label{eq:formula-u-2}
    u(z,T)
     = \int_{\T^n} g(x)\,d\sigma^z(x)
     + \int_{\T^n\times \R^n\times [0,T]} L(x,v) \,d\gamma^{z}(x,v,t).  
    \end{equation}
    Furthermore, $\gamma^z\in \cH(\sig^z,\delta_z;0,T)$, and
    \begin{equation}\label{eq:formula-minimizing-u}
    u(z,T) = \inf_{\gamma\in \cH(\delta_z;0,T)}\left( \int_{\T^n} g(x)\,d\sigma^\gamma(x)
     + \int_{\T^n\times \R^n\times [0,T]} L(x,v) \,d\gamma(x,v,t)\right).
    \end{equation}
    In other words, $\gamma^z$ is a minimizer of the above minimizing problem.
\end{thm}

\begin{rem}
 Theorem \ref{thm:rep2} is not new and a more general form (with a degenerate viscous term) already appeared in \cite[Theorem 1.1]{Gomes-Mitake-Tran}.
Nevertheless, the proof for the current first-order Hamilton--Jacobi equation is simpler and more direct, and we present it in this paper for completeness.   

Formula \eqref{eq:formula-minimizing-u} is a relaxed and generalized version of the classical optimal control formula.
See the discussion after the proof of  Theorem \ref{thm:rep2}.
For related problems, we refer the reader to \cite{Ishii-Mitake-Tran} and the references therein.

We note that Theorems \ref{thm:rep1}--\ref{thm:rep2} also hold in the whole space $\R^n$ -- that is, in the case where $H:\R^n\times \R^n\to \R$ and $g:\R^n \to \R$ are not assumed to be $\Z^n$-periodic -- under appropriate assumptions.
\end{rem}

Let us now consider the cell (ergodic) problem corresponding to \eqref{eq:HJ}.
For each $P\in \R^n$, the cell problem is
\begin{equation}\label{eq:cell-P}
    H(x,P+Dv_P(x)) = \ol H(P) \qquad \text{ in } \T^n.
\end{equation}
It was proved in \cite{LPV} that there exists a unique constant $\ol H(P)\in \R$ such that \eqref{eq:cell-P} has a viscosity solution $v_P\in C(\T^n)$.
Here, we only focus on the case where $P=0$, and
\begin{equation}\label{eq:cell}
    H(x,Dv_0(x)) = \ol H(0) \qquad \text{ in } \T^n.
\end{equation}
Let us recall the definition of Mather measures in the nonconvex setting in \cite[Theorem 1.3]{Cagnetti-Gomes-Tran}.

\begin{defn}[Mather measures in the nonconvex setting]
    Assume {\rm (A1)} and $P=0$.
    We say that a Radon probability measure $\mu$ on $\T^n\times \R^n$ is a Mather measure if the following items hold.
    \begin{itemize}
        \item[(a)] $\displaystyle \int_{\T^n\times \R^n} H(x,p)\,d\mu(x,p) = \ol H(0)=H(x,p)$ for $\mu$-a.e.~$(x,p)$.
        \item[(b)] $\displaystyle \int_{\T^n\times \R^n} D_pH(x,p)\cdot p\,d\mu(x,p) = 0$.
        \item[(c)] $\displaystyle \int_{\T^n\times \R^n} D_pH(x,p)\cdot D\phi(x)\,d\mu(x,p) = 0$ for all test functions $\phi\in C^2(\T^n)$.
    \end{itemize}
    
\end{defn}
If $H$ is convex, then the above definition coincides with the definition of Mather measures in the convex setting \cite{Mather, Mane, Evans-Gomes, Fathi, Tran-Yu} as proved in \cite[Theorem 1.3]{Cagnetti-Gomes-Tran}.
We use the ideas in Theorem \ref{thm:rep1} and large time averages to show the existence of Mather measures.

\begin{thm}\label{thm:Mather measures}
    Assume {\rm (A1)} and $P=0$.
    Then, the set of Mather measures is not empty.
\end{thm}

\begin{rem}
  Theorem \ref{thm:Mather measures} was already proved in \cite[Theorem 5.1]{Cagnetti-Gomes-Tran}.
  See also \cite{Guerra} for the space-time periodic nonconvex case.
Our approach here is a bit different and is based upon the ideas in Theorem \ref{thm:rep1} and large time averages.  

The weak KAM theory and the understanding of $\ol H$ in terms of dynamics in the nonconvex setting were listed as major open problems in \cite[Section 5]{Evans0}.
The results in \cite{Cagnetti-Gomes-Tran}, Theorem \ref{thm:Mather measures}, and Theorem \ref{thm:dissipative m} are the first steps in this direction.
Much needs to be studied, especially the properties of $\ol H$.
See \cite{Qian-Tran-Yu, tranbook} and the references therein for some decomposition formulas for $\ol H$ in some specific cases.

\end{rem}

Finally, we study the large time behavior of the solution to \eqref{eq:HJ} under (A1)--(A2).

\begin{thm}\label{thm:large time}
    Assume {\rm (A1)--(A2)} and $\ol H(0)=0$.
    Then, there exists a solution $v\in \Lip(\T^n)$ of \eqref{eq:cell} with $\ol H(0)=0$ such that
    \[
    \lim_{t\to \infty}\|u(\cdot,t)-v\|_{L^\infty(\T^n)}=0.
    \]
\end{thm}

\begin{rem}
    Large time behavior for \eqref{eq:HJ} was established for uniformly convex Hamiltonians in \cite{Fathi1, DS, II}; see also \cite{N-R} for some earlier results.
    For more general Hamiltonians including some possibly nonconvex ones, large time behavior for \eqref{eq:HJ} was obtained in \cite{BS1, BIM}.
    Roughly speaking, when $\ol H(0)=0$, \cite{BS1, BIM} only need to require some conditions similar to the strict convexity of $H$ near the $0$-sublevel set of $H$ (see conditions (A6)$_\pm$ in \cite{BIM} and condition (H4) in \cite{BS1} for further details).
    When the Hamiltonian is only convex, \cite{BS1} gives an example showing that large time behavior fails (see Example \ref{ex:H3} in Section \ref{sec:large time}).
    Intuitively, some strict convexity of $H$ near the $\ol H(0)$-sublevel set of $H$ was needed for the large time behavior results for general first-order Hamilton--Jacobi equations.
    
    For possibly degenerate viscous Hamilton--Jacobi equations with uniformly convex Hamiltonians, large time behavior was shown in \cite{CGMT}; see also \cite{Ley-Nguyen}.

    In Theorem \ref{thm:large time}, we deal with nonconvex Hamiltonians satisfying (A1)--(A2), which are different from the assumptions in \cite{BS1, BIM}.
    Note that (A2) is a bit similar to condition (H4)' in \cite{BS1}.
    Our result is new and complements the results in \cite{BS1, BIM}, which pushes further the study of large time behavior of nonconvex first-order Hamilton--Jacobi equations.
    See Examples \ref{ex:H1}--\ref{ex:H2} in Section \ref{sec:large time}.

    Our approach to prove Theorem \ref{thm:large time} is inspired by that of \cite{CGMT}.
    As we do not have the uniform convexity of $H$, we utilize the representation formula \eqref{eq:formula-u-1} and (A2) to get a one-sided control on $u^\ep_t$ and $u_t$, which is the most crucial step in the large time behavior proof.
\end{rem}

\subsection*{Organization of the paper}
The paper is organized as follows.
In Section \ref{sec:rep formulas}, we give the proofs of the representation formulas for the solutions, Theorem \ref{thm:rep1}--\ref{thm:rep2}.
We analyze Mather measures in the nonconvex setting and prove Theorem \ref{thm:Mather measures} in Section \ref{sec: MM}.
We also study dissipative measures and the invariance of Mather measures under the Hamiltonian flow in Theorem \ref{thm:dissipative m}.
The large time behavior result, Theorem \ref{thm:large time}, is proved in Section \ref{sec:large time}.
We also discuss further properties of the solution to \eqref{eq:HJ} and give some examples there.
Finally, Section \ref{sec:open problems} outlines open problems and questions for future research.

\section{Representation formulas} \label{sec:rep formulas}

\subsection{The general setting}
\begin{proof}[Proof of Theorem \ref{thm:rep1}]
    We compute
    \begin{align*}
        &\frac{d}{dt}\int_{\T^n} u^\ep \sig^{\ep,z}\,dx\\
        =\ &\int_{\T^n} \left(u^\ep_t \sig^{\ep,z}+u^\ep \sig^{\ep,z}_t\right)\,dx\\
        =\ &\int_{\T^n} \left((-H(x,Du^\ep)+\ep\Delta u^\ep) \sig^{\ep,z}+u^\ep (-\Div(D_pH\sig^{\ep,z})-\ep\Delta \sig^{\ep,z})\right)\,dx\\
        =\ &\int_{\T^n} \left(D_pH(x,Du^\ep)\cdot Du^\ep-H(x,Du^\ep)\right)\sig^{\ep,z}\,dx\\
        =\ &\int_{\T^n\times \R^n} \left(D_pH(x,p)\cdot p - H(x,p) \right) \,d\nu^{\ep,z}(x,p,t).
    \end{align*}
    Integrate this relation with respect to $t$ and note that $\sig^{\ep,z}(\cdot,T)=\delta_z$ and $u^\ep(\cdot,0)=g$ to yield
     \begin{align*}
    &u^\ep(z,T) \\
    = &\int_{\T^n} g(x)\sigma^{\ep,z}(x,0)\,dx + \int_0^T\int_{\T^n\times \R^n} \left(D_pH(x,p)\cdot p - H(x,p) \right)\, d\nu^{\ep,z}(x,p,t)dt\\
     = &\int_{\T^n} g(x)\sigma^{\ep,z}(x,0)\,dx + \int_{\T^n\times \R^n\times [0,T]} \left(D_pH(x,p)\cdot p - H(x,p) \right)\, d\mu^{\ep,z}(x,p,t).
    \end{align*}
    Recall that $\sigma^{\ep,z}(x,0)\,dx$ is a probability measure on $\T^n$, and $\nu^{\ep,z}(\cdot,\cdot,t)$ is a probability measure on $\T^n \times \ol B(0,C)$.
    Hence, $\mu^{\ep,z}$ is a nonnegative Radon measure, $\supp(\mu^{\ep,z}) \subset \T^n\times \ol B(0,C) \times [0,T]$, and $\mu^{\ep,z}(\T^n\times \ol B(0,C) \times [0,T])=T$.
    Thus, we can find a sequence $\{\ep_k\}\to 0$ such that $\sigma^{\ep_k,z}(x,0)\,dx$ converges to $d\sigma^z(x)$ and $\mu^{\ep_k,z}$ converges to $\mu^z$ weakly in the sense of measures. 
    It is clear that $\sig^z$ is a Radon probability measure on $\T^n$ and $\mu^z$ is a nonnegative Radon measure supported on $\T^n\times \ol B(0,C) \times [0,T]$ with $\mu^{z}(\T^n\times \ol B(0,C) \times [0,T])=T$.
    By letting $\ep=\ep_k \to 0$ in the above chain of equalities, we yield
       \begin{align*}
    u(z,T) 
     = \int_{\T^n} g(x)\,d\sigma^z(x) + \int_{\T^n\times \R^n\times [0,T]} \left(D_pH(x,p)\cdot p - H(x,p) \right) \,d\mu^{z}(x,p,t).
    \end{align*} 
\end{proof}

\begin{rem}\label{rem:char}
    The formula of $u(z,T)$ in the above theorem is a generalization of the formula from the characteristic method for Hamilton--Jacobi equations.
Let us recall quickly the method of characteristics.
For given $x_0\in \R^n$, we would like to compute the value of the solution $u$ along a curve emanating from $(x_0,0)$.
For $t\geq 0$, denote by
\begin{align*}
\begin{cases}
x(t): \text{ position of the curve at time } t, \\
 p(t): \text{ spatial gradient of } u \text{ at } (x(t), t),\\
  z(t): \text{ value of } u \text{ at } (x(t), t).
 \end{cases}
\end{align*}
The initial data is
\[
\begin{cases}
    x(0) = x_0, \\
    p(0) = Dg(x_0),\\
    z(0)=g(x_0).
\end{cases}
\]
For each such $x_0$, $(x(t),t)$ is called a characteristic.
We have the following system of ODE, which is often called the Hamiltonian system,
\[
\begin{cases}
    \dot {x}(t) = D_p H(x(t), p(t)), \\
    \dot {p}(t) = -D_x H(x(t), p(t)), \\
    \dot z(t) = p(t) \cdot D_p H(x(t), p(t)) - H(x(t), p(t)).
\end{cases}
\]

Once we can solve this system of ODE, then we can solve the Hamilton--Jacobi PDE \eqref{eq:HJ} at least locally.
The key point of the method of characteristics is that we can only define the solution if the characteristics $(x(t),t)$ do not cross each other.

If the characteristics do not cross yet at time $T$, and $x(T)=z$, then 
\[
u(z,T)=u(x(T),T)= g(x(0))+\int_0^T \left(p(t) \cdot D_p H(x(t), p(t)) - H(x(t), p(t))\right)\,dt.
\]
This formula is similar to \eqref{eq:formula-u-1}.

If the characteristics cross before time $T$, then we do not have the above formula.
It is then natural to see that the value of $u(z,T)$ is determined by a bunch of characteristics that all go through $(z,T)$.
All of this information is encoded in the measures $\sigma^z$ and $\mu^z$.
We emphasize that $\sigma^z$ and $\mu^z$ are dependent on the position $(z,T)$ in a highly nonlinear way.
\end{rem}

\subsection{The convex setting}

\begin{proof}[Proof of Theorem \ref{thm:rep2}]
We divide the proof into several steps.

\medskip

\noindent {\bf Step 1.}
    As $H$ is convex in $p$, we have
    \[
    D_pH(x,p)\cdot p - H(x,p) = L(x,D_pH(x,p)).
    \]
    Hence, thanks to \eqref{eq:formula-u-1},
        \begin{align*}
    u(z,T) 
     &= \int_{\T^n} g(x)\,d\sigma^z(x) + \int_{\T^n\times \R^n\times [0,T]} \left(D_pH(x,p)\cdot p - H(x,p) \right) \,d\mu^{z}(x,p,t)\\
     &=\int_{\T^n} g(x)\,d\sigma^z(x) + \int_{\T^n\times \R^n\times [0,T]} L(x,D_pH(x,p)) \,d\mu^{z}(x,p,t)\\
     &=\int_{\T^n} g(x)\,d\sigma^z(x) + \int_{\T^n\times \R^n\times [0,T]} L(x,v) \,d\gamma^{z}(x,v,t).
    \end{align*} 

\medskip

\noindent {\bf Step 2.}
We now prove that $\gamma^z \in \cH(\delta_z;0,T)$.
Fix $\varphi \in C^2(\T^n\times [0,T])$.
We compute
    \begin{align*}
        &\frac{d}{dt}\int_{\T^n} \varphi \sig^{\ep,z}\,dx\\
        =\ &\int_{\T^n} \left(\varphi_t \sig^{\ep,z}+\varphi \sig^{\ep,z}_t\right)\,dx\\
        =\ &\int_{\T^n} \left(\varphi_t \sig^{\ep,z}+\varphi (-\Div(D_pH\sig^{\ep,z})-\ep\Delta \sig^{\ep,z})\right)\,dx\\
        =\ &\int_{\T^n} \left(\varphi_t +D_pH(x,Du^\ep)\cdot D\varphi-\ep \Delta \varphi\right)\sig^{\ep,z}\,dx\\
        =\ &\int_{\T^n\times \R^n} \left(\varphi_t +D_pH(x,p)\cdot D\varphi  -\ep \Delta \varphi \right) \,d\nu^{\ep,z}(x,p,t).
    \end{align*}
    Integrate this relation with respect to $t$ to get
    \begin{multline*}
          \varphi(z,T) - \int_{\T^n} \varphi(x,0)\sigma^{\ep,z}(x,0)\,dx\\
    = \int_{\T^n\times \R^n\times [0,T]} \left(\varphi_t +D_pH(x,p)\cdot D\varphi  -\ep \Delta \varphi \right) \,d\mu^{\ep,z}(x,p,t).  
    \end{multline*}
    Let $\ep=\ep_k \to 0$ to imply
    \begin{align*}
    \varphi(z,T) - \int_{\T^n}\varphi(x,0)\,d\sigma^z(x)&=\int_{\T^n\times \R^n\times [0,T]} \left(\varphi_t +D_pH(x,p)\cdot D\varphi  \right) \,d\mu^{z}(x,p,t)\\
        &=\int_{\T^n\times \R^n\times [0,T]} (\varphi_t(x,t) + v\cdot D\varphi(x,t))\,d\gamma^z(x,v,t).
    \end{align*}
    Therefore, $\gamma^z \in \cH(\sigma^z,\delta_z;0,T) \subset \cH(\delta_z;0,T)$.

\medskip

\noindent {\bf Step 3.}
To conclude, we show that, for $\gamma\in \cH(\delta_z;0,T)$,
    \[
    u(z,T) \leq  \int_{\T^n} g(x)\,d\sigma^\gamma(x)
     + \int_{\T^n\times \R^n\times [0,T]} L(x,v) \,d\gamma(x,v,t).
    \]
Let $\xi$ be a standard mollifier in $\R^{n+1}$, that is, 
\[
\xi \in C_c^\infty(\R^{n+1},[0,\infty)), \quad \supp(\xi) \subset B(0,1), \quad \int_{\R^{n+1}} \xi(z)\,dz=1.
\]
For $\alpha\in (0,1)$, let $\xi^\alpha(z) = \alpha^{-(n+1)} \xi(z/\alpha)$ for $z\in \R^{n+1}$.
For $(x,t)\in \T^n\times [0,\infty)$, we define
\[
\tilde u^\alpha(x,t) = (\xi^\alpha\star u)(x,t+\alpha)=\int_{\R^{n+1}} \xi^\alpha(y,s) u(x-y,t+\alpha-s)\,dyds.
\]
As $u$ is Lipschitz with constant $C$, we have that $\tilde u^\alpha$ is Lipschitz with constant $C$.
Further, $\tilde u^\alpha$ is smooth and $\tilde u^\alpha \to u$ locally uniformly as $\alpha\to 0$.
By \eqref{eq:HJ},
\begin{align*}
    0&= \tilde u^\alpha_t(x,t)+\int_{B(0,\alpha)} \xi^\alpha(y,s)H(x-y,D u(x-y,t+\alpha-s))\,dyds\\
    &\geq \tilde u^\alpha_t(x,t)+\int_{B(0,\alpha)} \xi^\alpha(y,s)H(x,D u(x-y,t+\alpha-s))\,dyds - C\alpha\\
    &\geq \tilde u^\alpha_t(x,t)+H\left(x,\int_{B(0,\alpha)} \xi^\alpha(y,s)D u(x-y,t+\alpha-s)\,dyds\right) - C\alpha\\
    &= \tilde u^\alpha_t(x,t)+H(x, D \tilde u^\alpha(x,t)) - C\alpha.
\end{align*}
We used Jensen's inequality in the second last line above.
Hence, for $v\in \R^n$,
\begin{align*}
    C\alpha&\geq \tilde u^\alpha_t(x,t)+H(x, D \tilde u^\alpha(x,t))\\
    &\geq \tilde u^\alpha_t(x,t)+v\cdot D \tilde u^\alpha(x,t)-L(x,v).
\end{align*}
Integrate this inequality with respect to $d\gam(x,v,t)$ and use the fact that $\gamma\in \cH(\sigma^\gamma,\delta_z;0,T)$ to yield
\[
\tilde u^\alpha(z,T) - \int_{\T^n}\tilde u^\alpha(x,0)\,d\sigma^\gamma(x) - \int_{\T^n \times \R^n\times [0,T]} L(x,v)\,d\gamma(x,v,t) \leq CT\alpha.
\]
We let $\alpha \to 0$ to conclude that
    \[
    u(z,T) \leq  \int_{\T^n} g(x)\,d\sigma^\gamma(x)
     + \int_{\T^n\times \R^n\times [0,T]} L(x,v) \,d\gamma(x,v,t).
    \]
\end{proof}

We next make some connections between Theorem \ref{thm:rep2} and the optimal control formula.
Since we are in the convex setting, we have
\begin{equation}\label{eq:OC}
u(z,T) = \inf \left\{\int_0^T L(\zeta(t),\dot\zeta(t))\,dt + g(\zeta(0))\,:\, \zeta\in \AC([0,T],\T^n), \zeta(T)=z\right\}.
\end{equation}
Assume that $\xi$ is a minimizer to \eqref{eq:OC}, that is, $\xi(T)=z$, and
\[
u(z,T) = \int_0^T L(\xi(t),\dot\xi(t))\,dt + g(\xi(0)).
\]
Under some appropriated conditions, we have that $\xi \in C^2([0,T])$ and it satisfies the Euler--Lagrange equations
\[
\frac{d}{dt}\left(D_v L(\xi(t),\dot \xi(t)) \right) = D_x L(\xi(t),\dot \xi(t)) \quad \text{ for } t\in (0,T).
\]
We define the  measure $\gamma^\xi$ corresponding to $\xi$ as
\[
d\gamma^\xi(x,v,t) = \delta_{(\xi(t),\dot \xi(t))}dt.
\]
We can see that $\gamma^\xi \in \cH(\delta_{\xi(0)},\delta_z;0,T)$ as
\begin{align*}
    &\int_{\T^n\times \R^n\times [0,T]} (\varphi_t(x,t) + v\cdot D\varphi(x,t))\,d\gamma^\xi(x,v,t)\\
    =\ & \int_0^T \left( \varphi_t(\xi(t),t) + \dot \xi(t)\cdot D\varphi(\xi(t),t)\right)\,dt
    =\varphi(z,T) - \varphi(\xi(0),0).
\end{align*}
for all test functions $\varphi \in C^2(\T^n\times [0,T])$.
Thus, $\gam^\xi$ is a minimizer to \eqref{eq:formula-minimizing-u}.
This shows that each minimizer of \eqref{eq:OC} corresponds to a minimizer of \eqref{eq:formula-minimizing-u}, a more relaxed problem.
The converse is not true in general as if $\gamma_1, \gamma_2$ are two minimizers of \eqref{eq:formula-minimizing-u}, then convex combinations of $\gamma_1, \gamma_2$, that is, $s\gamma_1+(1-s)\gamma_2$ for $s\in [0,1]$, are also minimizers.
This leads naturally to a selection problem: Do we have that $\{\mu^{\ep,z}\}$ converges weakly in the sense of measures to a unique limit as $\ep \to 0$?
See Questions \ref{quest1}--\ref{quest2} in Section \ref{sec:open problems}.

\section{Mather measures in the nonconvex setting} \label{sec: MM}

\subsection{Existence of Mather measures}
We first approximate the cell problem \eqref{eq:cell} by the vanishing viscosity process.
For $\ep\in (0,1)$, we consider
\begin{equation}\label{eq:cell-ep}
    H(x,Dv^\ep) = \ol H^\ep(0) + \ep \Delta v^\ep \qquad \text{ in } \T^n.
\end{equation}
There exists a unique constant $\ol H^\ep(0)\in \R$ such that \eqref{eq:cell-ep} has a solution $v^\ep \in C^2(\T^n)$.
Further, $v^\ep$ is the unique solution to \eqref{eq:cell-ep} up to additive constants.
We normalize it so that $v^\ep(0)=0$.
Equation \eqref{eq:cell-ep} was studied in the convex setting in \cite{Gomes, IS}.
Let us summarize some basic properties of \eqref{eq:cell-ep} and \eqref{eq:cell} in the following proposition (see \cite[Theorem 2.1]{Cagnetti-Gomes-Tran} or \cite[Proposition 5.5]{Le-Mitake-Tran} for example).

\begin{prop}\label{prop:basic}
    Assume {\rm (A1)}.
    For $\ep\in (0,1)$, let $v^\ep$ be the solution to \eqref{eq:cell-ep} with $v^\ep(0)=0$.
    The following properties hold.
    \begin{itemize}
        \item[(i)] There exists $C>0$ independent of $\ep\in (0,1)$ such that
        \[
        \|Dv^\ep\|_{L^\infty(\T^n)} \leq C.
        \]
        \item[(ii)] There exists a sequence $\{\ep_k\} \to 0$ such that $v^{\ep_k} \to v$ uniformly on $\T^n$, and $v\in \Lip(\T^n)$ is a solution to \eqref{eq:cell}.

        \item[(iii)] There exists $C>0$ independent of $\ep\in (0,1)$ such that
        \[
        \left| \ol H^\ep(0) - \ol H(0) \right| \leq C \ep^{1/2}.
        \]
    \end{itemize}
\end{prop}
In the general setting, the bound $\left| \overline{H}^\varepsilon(0) - \overline{H}(0) \right| \leq C \varepsilon^{1/2}$ is the best we can expect.
When $H$ is uniformly convex in $p$, the convergence rate improves to $C\varepsilon$ (see \cite{Evans, AIPS, Yu} for the classical mechanic Hamiltonian, and \cite{Tu-Zhang} for general uniformly convex Hamiltonians).

\medskip

We note that $u^\ep(x,t)=v^\ep(x) - \ol H^\ep(0)t$ is the solution to \eqref{eq:HJ-ep} with initial data $g=v^\ep$, that is,
\begin{equation*}
    \begin{cases}
        u^\ep_t + H(x,Du^\ep)=\ep \Delta u^\ep \qquad &\text{ in } \T^n \times (0,\infty),\\
        u^\ep(x,0)=v^\ep(x) \qquad &\text{ on } \T^n.
    \end{cases}
\end{equation*}
And $u^{\ep_k} \to u$ locally uniformly on $\T^n\times [0,\infty)$ as $\ep_k \to 0$, where $u(x,t) = v(x) - \ol H(0) t$ solves
\begin{equation*}
    \begin{cases}
        u_t + H(x,Du)=0 \qquad &\text{ in } \T^n \times (0,\infty),\\
        u(x,0)=v(x) \qquad &\text{ on } \T^n.
    \end{cases}
\end{equation*}

As $\|Du^\ep\|_{L^\infty(\T^n\times [0,\infty))}=\|Dv^\ep\|_{L^\infty(\T^n)} \leq C$, for each $t\in [0,T]$, let $\nu^{\ep,z}(\cdot,\cdot,t)$ be the probability measure on $\T^n \times \ol B(0,C)$ such that, for all $\psi \in C(\T^n\times \R^n)$,
\begin{align*}
\int_{\T^n\times \R^n}\psi(x,p)\,d\nu^{\ep,z}(x,p,t) 
&= \int_{\T^n}\psi(x,Du^\ep(x,t))\sigma^{\ep,z}(x,t)\,dx\\
&= \int_{\T^n}\psi(x,Dv^\ep(x))\sigma^{\ep,z}(x,t)\,dx.
\end{align*}
Denote by $\mu^{\ep,z,T}$ the measure on $\T^n\times \ol B(0,C)\times [0,T]$ such that 
\[
d\mu^{\ep,z,T}(x,p,t)=d\nu^{\ep,z}(x,p,t)dt.
\]

\begin{lem}\label{lem:NAM}
    Assume {\rm (A1)}.
    There exists $C>0$ independent of $\ep\in (0,1)$ and $T>0$ such that
    \[
    \ep \int_0^T \int_{\T^n} |D^2u^\ep|^2 \sig^{\ep,z}\,dxdt=\ep \int_0^T \int_{\T^n} |D^2v^\ep|^2 \sig^{\ep,z}\,dxdt \leq C(1+T).
    \]
\end{lem}

\begin{proof}
    Let $\phi=|Du^\ep|^2/2$.
    Then $\phi$ satisfies
    \[
    \phi_t + D_p H(x,Du^\ep)\cdot D\phi - \ep \Delta \phi + \ep |D^2 u^\ep|^2 + D_xH(x,Du^\ep)\cdot Du^\ep=0.
    \]
    Hence,
    \[
    \cL^\ep[\phi] + \ep  |D^2 u^\ep|^2 \leq \left| D_xH(x,Du^\ep)\cdot Du^\ep \right| \leq C.
    \]
    Multiply this inequality by $\sig^{\ep,z}$ and integrate to imply
    \[
     \ep \int_0^T \int_{\T^n} |D^2u^\ep|^2 \sig^{\ep,z}\,dxdt \leq CT + \int_{\T^n}\phi(x,0)\sig^{\ep,z}(x,0)\,dx - \phi(z,T) \leq C(1+T).
    \]
\end{proof}

By passing to a subsequence if needed, we assume that $\mu^{\ep_k,z,T}$ converges weakly in the sense of measures to $\mu^{z,T} \in \mathcal{R}^+(\T^n\times \ol B(0,C)\times [0,T])$.

\begin{lem}\label{lem:energy level}
    Assume {\rm (A1)}.
    We have
    \[
    \int_{\T^n\times \R^n\times [0,T]} \left| H(x,p) - \ol H(0)\right|^2\,d\mu^{z,T}(x,p,t)=0.
    \]
\end{lem}

\begin{proof}
    By \eqref{eq:cell-ep},
    \[
    \left|H(x,Dv^\ep) - \ol H^\ep(0)\right|^2 = \ep^2 \left(\Delta v^\ep\right)^2.
    \]
    Multiply this by $\sig^{\ep,z}$ and integrate to yield
    \begin{align*}
        &\int_{\T^n\times \R^n\times [0,T]} \left| H(x,p) - \ol H^\ep(0)\right|^2\,d\mu^{\ep,z,T}(x,p,t)\\
       =\ & \ep^2 \int_0^T \int_{\T^n} |\Delta v^\ep|^2 \sig^{\ep,z}\,dxdt
       \leq  n\ep^2 \int_0^T \int_{\T^n} |D^2 v^\ep|^2 \sig^{\ep,z}\,dxdt \leq C\ep(1+T).
    \end{align*}
    Letting $\ep=\ep_k \to 0$ and using the fact that $\left| \ol H^\ep(0) - \ol H(0) \right| \leq C \ep^{1/2}$, we conclude
       \[
    \int_{\T^n\times \R^n\times [0,T]} \left| H(x,p) - \ol H(0)\right|^2\,d\mu^{z,T}(x,p,t)=0.
    \]
\end{proof}

Recall that $\mu^{z,T} \in \mathcal{R}^+(\T^n\times \ol B(0,C)\times [0,T])$ and $\mu^{z,T}(\T^n\times \ol B(0,C)\times [0,T])=T$.
For any Borel measurable set $A\subset \T^n\times \ol B(0,C)$, denote by
\[
\tilde \mu^{z,T}(A) = \frac{1}{T}\mu^{z,T}(A\times [0,T]).
\]
Basically, we take the large time average of $\mu^{z,T}$ to get $\tilde \mu^{z,T}\in \cP(\T^n\times \ol B(0,C))$.
Naturally, as $T\to\infty$, we should expect ergodicity, that is, the appearance of Mather measures.

\begin{thm}\label{thm:M-measure}
    Take a sequence $T_k \to \infty$ such that $\tilde \mu^{z,T_k}$ converges weakly in the sense of measures to $\tilde \mu\in \cP(\T^n\times \ol B(0,C))$.
    Then, $\tilde \mu$ is a Mather measure.
\end{thm}
\begin{proof}
    Firstly, by Lemma \ref{lem:energy level},
     \[
    \int_{\T^n\times \R^n} \left| H(x,p) - \ol H(0)\right|^2\,d\tilde\mu^{z,T}(x,p)=0.
    \]
    Let $T=T_k \to \infty$ to yield that
     \begin{equation}\label{eq:MM1}
    \int_{\T^n\times \R^n} \left| H(x,p) - \ol H(0)\right|^2\,d\tilde\mu(x,p)=0.
    \end{equation}
    In particular,
    \[
    \int_{\T^n\times \R^n}  H(x,p) \,d\tilde\mu(x,p)= \ol H(0)= H(x,p) \quad \text{ for $\tilde \mu$-a.e. } (x,p).
    \]

    Next, in light of Theorem \ref{thm:rep1},
    \begin{align*}
    u(z,T) 
     = \int_{\T^n} v(x)\,d\sigma^{z,T}(x) + \int_{\T^n\times \R^n\times [0,T]} \left(D_pH(x,p)\cdot p - H(x,p) \right) \,d\mu^{z,T}(x,p,t).
    \end{align*} 
    As $u(z,T)=v(z) - \ol H(0)T$, we divide both sides by $T$ and rearrange to get
    \[
    \frac{1}{T}v(z) -\frac{1}{T} \int_{\T^n} v(x)\,d\sigma^{z,T}(x) =
    \ol H(0)+\int_{\T^n\times \R^n} \left(D_pH(x,p)\cdot p - H(x,p) \right) \,d\tilde \mu^{z,T}(x,p).
    \]
    Thanks to \eqref{eq:MM1},
    \[
    \frac{1}{T}v(z) -\frac{1}{T} \int_{\T^n} v(x)\,d\sigma^{z,T}(x) =
    \int_{\T^n\times \R^n} D_pH(x,p)\cdot p  \,d\tilde \mu^{z,T}(x,p).
    \]
    By letting $T=T_k\to\infty$, we deduce
    \begin{equation}\label{eq:MM2}
     \int_{\T^n\times \R^n}   D_pH(x,p)\cdot p  \,d\tilde \mu(x,p)=0.
    \end{equation}

    Finally, fix $\phi\in C^2(\T^n)$.
    Multiply \eqref{eq:sigma-ep} by $\phi$ and integrate to imply
    \[
    \int_{\T^n\times \R^n\times [0,T]} \left( D_pH(x,p)\cdot D\phi(x) -\ep \Delta \phi \right)\,d\mu^{\ep,z,T}(x,p,t) = \phi(z) - \int_{\T^n} \phi \sig^{\ep,z}(x,0)\,dx.
    \]
    Divide both sides by $T$ and let $\ep=\ep_k \to 0$, $T=T_k \to \infty$ in this order to conclude that
      \begin{equation}\label{eq:MM3}
     \int_{\T^n\times \R^n}   D_pH(x,p)\cdot D\phi(x)  \,d\tilde \mu(x,p)=0.
    \end{equation}
    By \eqref{eq:MM1}--\eqref{eq:MM3}, we have that $\tilde \mu$ is a Mather measure.
    The proof is complete.
    \end{proof}

    \begin{proof}[Proof of Theorem \ref{thm:Mather measures}]
    We see that Theorem \ref{thm:M-measure} implies Theorem \ref{thm:Mather measures} immediately.
   \end{proof}

    Let us discuss the case where we know a bit more about the Mather measures.
    
    \begin{lem}\label{lem:MM ex}
    Assume {\rm (A1)--(A2)} and $\ol H(0)=0$.
    Let $\mu$ be a Mather measure.
    Then,
    \[
    \supp(\mu) \subset \left\{(x,p)\in \T^n\times \R^n\,:\,H(x,p)=D_pH(x,p)\cdot p=0 \right\}.
    \]
    
    \end{lem}

    \begin{proof}
    By property (a) of Mather measures, for $\mu$-a.e. $(x,p)$,
    \[
    H(x,p)=\int_{\T^n\times \R^n} H(x,p)\,d\mu(x,p) = \ol H(0)=0,
    \]
    which, together with (A2), yields
    \[
    D_pH(x,p)\cdot p \geq (\theta+1)H(x,p)=0.
    \]
    Therefore,
    \[
    \int_{\T^n\times \R^n} D_pH(x,p)\cdot p\,d\mu(x,p) \geq 0.
    \]
    By property (b) of Mather measures, we must have equality in the above.
    Hence, for $\mu$-a.e. $(x,p)$,
    \[
    H(x,p)=D_pH(x,p)\cdot p=0.
    \]
    The proof is complete.
    \end{proof}

    \begin{ex}\label{ex:H0}
    Assume
    \[
    H(x,p)=H(p)=|p|^4-|p|^2 \qquad \text{ for } (x,p) \in \T^n\times \R^n.
    \]
    Then, $\ol H(0)=0$ as $H(0)=0$.
    It is clear that {\rm (A1)} holds.
    We compute
    \begin{align*}
        &D_pH(x,p)\cdot p - H(x,p)=D_pH(p)\cdot p - H(p)\\
        =\ &(4|p|^2p - 2p)\cdot p - (|p|^4-|p|^2)=3|p|^4-|p|^2 \geq H(p).
    \end{align*}
    We have {\rm (A2)} holds with $\theta=1$.
    In this situation,
    \[
    H(p)=D_pH(p)\cdot p = 0 \quad \iff \quad p=0.
    \]
    By Lemma \ref{lem:MM ex}, we have, if $\mu$ is a Mather measure, then
    \[
    \supp(\mu) \subset \T^n \times \{0\}.
    \]
    On the other hand, for any probability measure $\mu$ of the form
    \[
    \mu(x,p)=\nu(x)\delta_0(p),
    \]
    where $\nu\in\cP(\T^n)$, we see that $\mu$ is a Mather measure.
    Thus, we have the complete description of Mather measures in this simple case.
    \end{ex}

    \subsection{Dissipative measures}
     We study whether the Mather measures are invariant under the Hamiltonian flow.
    In fact, in the general nonconvex setting, this invariant property might fail, and a dissipation arises, which is described by a positive semidefinite matrix of Borel measures in \cite[Theorem 5.1]{Cagnetti-Gomes-Tran}.
    We give a different proof of \cite[Theorem 5.1]{Cagnetti-Gomes-Tran} in this section based on the tools developed in this paper.

\medskip
    We start with some definitions.

    \begin{defn}[Poisson's bracket]
        Let $F,G\in C^1(\T^n\times \R^n)$.
        The Poisson bracket between $F$ and $G$ is defined as
        \[
        \{F,G\} = D_pF\cdot D_xG - D_xF \cdot D_pG.
        \]
    \end{defn}

    \begin{defn}[Invariance under the Hamiltonian flow]
Let $\mu$ be a Mather measure.
We say that $\mu$ is invariant under the Hamiltonian flow if for any $\psi\in C^2_c(\T^n\times \R^n)$,
\[
\int_{\T^n\times \R^n}\{H,\psi\}(x,p)\,d\mu(x,p)=0.
\]
    \end{defn}

    \begin{thm}\label{thm:dissipative m}
        Let $\tilde \mu$ be the Mather measure constructed in Theorem \ref{thm:M-measure}.
        Then, there exists a nonnegative symmetric matrix $(m_{ij})_{1\leq i,j \leq n}$ of Borel measures on $\T^n \times \R^n$ such that, for any $\psi\in C^2_c(\T^n\times \R^n)$,
\[
\int_{\T^n\times \R^n}\{H,\psi\}(x,p)\,d\tilde \mu(x,p)=\int_{\T^n\times \R^n}\psi_{p_i p_j}(x,p)\, dm_{ij}(x,p).
\]
We call $(m_{ij})_{1\leq i,j \leq n}$ the matrix of dissipative measures.
    \end{thm}
    
We note that Theorem \ref{thm:Mather measures} and Theorem \ref{thm:dissipative m} combined together is exactly \cite[Theorem 5.1]{Cagnetti-Gomes-Tran}.

    \begin{proof}
    Recall that $u^\ep(x,t)=v^\ep(x) - \ol H^\ep(0)t$ is the solution to \eqref{eq:HJ-ep} with initial data $g=v^\ep$.
    In particular, $u^\ep_t=-\ol H^\ep(0)$, $Du^\ep=Dv^\ep$, and $\Delta u^\ep =\Delta v^\ep$.
    Differentiate \eqref{eq:HJ-ep} with respect to $x_i$ to get
    \[
    H_{x_i}(x,Dv^\ep) + H_{p_j}(x,Dv^\ep) v^\ep_{x_i x_j} =\ep \Delta v^\ep_{x_i}.
    \]
    Multiply this relation by $\psi_{p_i}(x,Dv^\ep) \sig^{\ep,z}$ and integrate
    \begin{multline}\label{eq:diss-1}
        \int_0^T\int_{\T^n} D_xH\cdot D_p\psi \sig^{\ep,z}\,dxdt +  \int_0^T\int_{\T^n} H_{p_j}\psi_{p_i} v^\ep_{x_i x_j} \sig^{\ep,z}\,dxdt\\
        = \ep  \int_0^T\int_{\T^n}
        \psi_{p_i}\Delta v^\ep_{x_i} \sig^{\ep,z}\,dxdt.
    \end{multline}
    Next, multiply the nonlinear adjoint equation \eqref{eq:sigma-ep} by $\psi(x,Dv^\ep)$ and integrate by parts to imply
    \begin{multline*}
        \int_{\T^n} \psi(x,Dv^\ep)\sig^{\ep,z}(x,0)\,dx - \psi(z,Dv^\ep(z))+ \int_0^T\int_{\T^n} D_pH\cdot D_x(\psi(x,Dv^\ep)) \sig^{\ep,z}\,dxdt\\
        =\ep  \int_0^T\int_{\T^n} \Delta_x(\psi(x,Dv^\ep)) \sig^{\ep,z}\,dxdt.
    \end{multline*}
    Hence,
    \begin{multline}\label{eq:diss-2}
        \int_{\T^n} \psi(x,Dv^\ep)\sig^{\ep,z}(x,0)\,dx - \psi(z,Dv^\ep(z))+ \int_0^T\int_{\T^n} D_pH\cdot D_x\psi(x,Dv^\ep) \sig^{\ep,z}\,dxdt\\
        +\int_0^T\int_{\T^n} H_{p_j} \psi_{p_i} v^\ep_{x_i x_j} \sig^{\ep,z}\,dxdt
        =\ep  \int_0^T\int_{\T^n} \psi_{p_i p_j}(x,Dv^\ep)v^{\ep}_{x_i x_k} v^\ep_{x_j x_k}\sig^{\ep,z}\,dxdt\\
        +\ep  \int_0^T\int_{\T^n} \left(\psi_{x_i x_i} +2\psi_{x_i p_j}v^\ep_{x_i x_j} + \psi_{p_i}\Delta v^\ep_{x_i}\right)\sig^{\ep,z}\,dxdt.
    \end{multline}
    Take the difference of \eqref{eq:diss-2} and \eqref{eq:diss-1} to get
     \begin{multline*}
        \int_0^T\int_{\T^n} \{H,\psi\}(x,Dv^\ep) \sig^{\ep,z}\,dxdt
        =  \ep  \int_0^T\int_{\T^n} \psi_{p_i p_j}(x,Dv^\ep)v^{\ep}_{x_i x_k} v^\ep_{x_j x_k}\sig^{\ep,z}\,dxdt\\
        +\psi(z,Dv^\ep(z))
        -\int_{\T^n} \psi\sig^{\ep,z}(x,0)\,dx 
        +\ep  \int_0^T\int_{\T^n} \left(\psi_{x_i x_i} +2\psi_{x_i p_j}v^\ep_{x_i x_j} \right)\sig^{\ep,z}\,dxdt.
    \end{multline*}
    As $\|Dv^\ep\|_{L^\infty(\T^n)} \leq C$,
    \[
    \left|\psi(z,Dv^\ep(z))
        -\int_{\T^n} \psi\sig^{\ep,z}(x,0)\,dx \right| \leq C,
    \]
    and
    \begin{align*}
        &\ep  \int_0^T\int_{\T^n} \left(|\psi_{x_i x_i}| +2|\psi_{x_i p_j}v^\ep_{x_i x_j}| \right)\sig^{\ep,z}\,dxdt\\
        \leq \ &C\ep  \int_0^T\int_{\T^n} \left(1+|D^2 v^\ep| \right)\sig^{\ep,z}\,dxdt
        \leq CT (\ep +\ep^{1/2}) \leq CT\ep^{1/2}.
    \end{align*}
    We used Lemma \ref{lem:NAM} in the second last inequality above.
    Therefore,
    \begin{multline}\label{eq:diss-3}
       \frac{1}{T} \int_{\T^n\times \R^n\times [0,T]} \{H,\psi\}(x,p) \,d\mu^{\ep,z,T}(x,p,t)\\
        = \frac{\ep}{T}  \int_{\T^n\times \R^n\times [0,T]} \psi_{p_i p_j}(x,p)v^{\ep}_{x_i x_k} v^\ep_{x_j x_k}\,d\mu^{\ep,z,T}(x,p,t)
        +O\left(\frac{1}{T}\right) + O\left(\ep^{1/2}\right).
    \end{multline}
    We note that the first term on the right-hand side above does not vanish as $\ep \to 0$ in general.
    In fact, it is only bounded thanks to Lemma \ref{lem:NAM}.
    By taking a subsequence of $\{\ep_k\}$ if needed, we can assume that, for every $\phi\in C_c(\T^n\times \R^n)$,
    \[
    \lim_{k\to \infty}\frac{\ep_k}{T}  \int_{\T^n\times \R^n\times [0,T]} \phi(x,p)v^{\ep_k}_{x_i x_k} v^\ep_{x_j x_k}\,d\mu^{\ep_k,z,T}(x,p,t)
    =\int_{\T^n\times \R^n} \phi(x,p) \,dm^{z,T}_{ij}(x,p),
    \]
    where $(m^{z,T}_{ij})_{1\leq i,j \leq n}$ is a nonnegative symmetric matrix of Borel measures on $\T^n \times \ol B(0,C)$.
    In light of Lemma \ref{lem:NAM}, $(m^{z,T}_{ij})_{1\leq i,j \leq n}$ is bounded.

    \medskip
    
    We let $\ep=\ep_k \to 0$ in \eqref{eq:diss-3} to deduce
    \[
    \int_{\T^n\times \R^n} \{H,\psi\}(x,p) \,d\tilde \mu^{z,T}(x,p)
    =\int_{\T^n\times \R^n} \psi_{p_i p_j}(x,p) \,dm^{z,T}_{ij}(x,p)+O\left(\frac{1}{T}\right).
    \]
    Finally, let $T=T_k \to \infty$ and pass to a subsequence if necessary to get
    \[
    \int_{\T^n\times \R^n} \{H,\psi\}(x,p) \,d\tilde \mu(x,p)
    =\int_{\T^n\times \R^n} \psi_{p_i p_j}(x,p) \,dm_{ij}(x,p).
    \]
    Here, $(m^{z,T_k}_{ij})_{1\leq i,j \leq n}$ converges weakly in the sense of measures  to $(m_{ij})_{1\leq i,j \leq n}$ as $T_k \to \infty$.
    Of course, $(m_{ij})_{1\leq i,j \leq n}$ is a nonnegative symmetric matrix of Borel measures on $\T^n \times \ol B(0,C)$, and $(m_{ij})_{1\leq i,j \leq n}$ is bounded.
    \end{proof}

    From Theorem \ref{thm:dissipative m}, we see that the Mather measure $\tilde \mu$ is invariant under the Hamiltonian flow if and only if the dissipative measures $(m_{ij})_{1\leq i,j \leq n}$ vanish.
    As noted in \cite{Cagnetti-Gomes-Tran}, this is not always the case, and in general, $(m_{ij})_{1\leq i,j \leq n}$ records the gradient jumps of the solution to \eqref{eq:HJ} along the shock curves.

    Let us now give a class of Hamiltonians where the dissipative measures $(m_{ij})_{1\leq i,j \leq n}$ vanish, which is different from the various classes given in \cite{Cagnetti-Gomes-Tran}.

    \begin{lem}\label{lem:new diss}
        Assume that $H(x,p)=c(x)K(p)$ for $(x,p)\in \T^n\times \R^n$, where $c\in C^2(\T^n,(0,\infty))$ and $K\in C^2(\R^n)$, satisfies {\rm (A1)}.
        Assume further that $\ol H(0)=0$.
        Let $\tilde \mu$ be the Mather measure constructed in Theorem \ref{thm:M-measure}.
        Then, $(m_{ij})_{1\leq i,j \leq n}$ vanish.
        In other words, $\tilde \mu$ is invariant under the Hamiltonian flow.
        
    \end{lem}

    \begin{proof}
        As $H(x,p)=c(x)K(p)$, $D_x H(x,p)=K(p)Dc(x)$.
        By the definition of Mather measures and the fact that $c\in C^2(\T^n,(0,\infty))$, we have, for $\tilde \mu$-a.e. $(x,p)$,
        \[
        H(x,p)=K(p)=0.
        \]
        Hence, 
        for $\tilde \mu$-a.e. $(x,p)$,
        \begin{equation}\label{eq:new diss-1}
        D_xH(x,p) = K(p)Dc(x)=0.
        \end{equation}
        Take $\psi(x,p)=|p|^2/2$ for $(x,p)\in \T^n\times \R^n$.
        Then, thanks to \eqref{eq:new diss-1},
        \[
        \int_{\T^n\times \R^n} \{H,\psi\}(x,p) \,d\tilde \mu(x,p)
    =\int_{\T^n\times \R^n} -D_xH(x,p)\cdot p \,d\tilde \mu(x,p)
    =0.
        \]
        Combine this with Theorem \ref{thm:dissipative m} to yield
        \[
        \int_{\T^n \times \R^n} dm_{ii}(x,p)=0.
        \]
        As $(m_{ij})_{1\leq i,j \leq n}$ is a nonnegative symmetric matrix of Borel measures on $\T^n \times \ol B(0,C)$, we conclude that $(m_{ij})_{1\leq i,j \leq n}=0$.
    \end{proof}

    \begin{rem}
        Some comments about Lemma \ref{lem:new diss} are in order.

        \smallskip
       
        Firstly, we give an example of $H$ satisfying the required assumptions.
            Let $H(x,p)=c(x)K(p)$ for $(x,p)\in \T^n\times \R^n$, where $c\in C^2(\T^n,(0,\infty))$ and $K\in C^2(\R^n)$.
            Assume further that $K(0)=0$ and $K$ is superlinear, that is,
            \[
            \lim_{|p|\to\infty} \frac{K(p)}{|p|}=+\infty.
            \]
        Then, (A1) holds as
        \begin{align*}
        &\lim_{|p|\to \infty} \min_{x\in \T^n} \left(\frac{1}{2}H(x,p)^2 + D_xH(x,p)\cdot p \right)\\
        =\ &\lim_{|p|\to \infty} \left[K(p)^2\times\min_{x\in \T^n}  \left(\frac{1}{2}c(x)^2 + Dc(x)\cdot \frac{p}{K(p)}\right)\right]
        =+\infty.
        \end{align*}
        Besides, it is clear that $\ol H(0)=0$ as $v_0\equiv 0$ is a solution to \eqref{eq:cell}.    
        We note that $K$ is not required to be convex and its zero level set $\{p\in \R^n\,:\, K(p)=0\}$ can be extremely complicated.
        Nevertheless, the dissipative measures vanish in this case.

        \smallskip

        Secondly, through the proof of Lemma \ref{lem:new diss}, we see that $(m_{ij})_{1\leq i,j \leq n}=0$ if and only if
        \[
        \int_{\T^n\times \R^n} D_xH(x,p)\cdot p \,d\tilde \mu(x,p)=0.
        \]
        
    \end{rem}

\section{Large time behavior in a nonconvex setting}\label{sec:large time}
\subsection{Large time behavior}
In this section, we always assume (A1)--(A2) and $\ol H(0)=0$.
Let $v_0 \in \Lip(\T^n)$ be a solution to \eqref{eq:cell}.
Then, $v_0+C$ is a stationary solution to \eqref{eq:HJ} for any $C\in \R$.
By the usual comparison principle, we have
\begin{equation}\label{eq:LT1}
    v_0-\|v_0\|_{L^\infty(\T^n)}-\|g\|_{L^\infty(\T^n)} \leq u \leq v_0+\|v_0\|_{L^\infty(\T^n)}+\|g\|_{L^\infty(\T^n)}.
\end{equation}

We first use large time average to give a representation formula for $u^\ep_t$.

\begin{lem}\label{lem:LT1}
    Assume {\rm (A1)--(A2)}.
    Then,
    \[
    u^\ep_t(z,T) = \frac{1}{T}\int_0^T\int_{\T^n}\left(-H(x,Du^\ep)+\ep \Delta u^\ep\right)\sig^{\ep,z}\,dxdt.
    \]
\end{lem}

\begin{proof}
    Differentiate \eqref{eq:HJ-ep} with respect to $t$ to get
    \[
    \cL^\ep[u^\ep_t] = u^\ep_{tt}+D_pH(x,Du^\ep)\cdot Du^\ep_t - \ep \Delta u^\ep_t = 0.
    \]
    We multiply this relation by $\sig^{\ep,z}$ and integrate over $\T^n$
    \[
    \frac{d}{dt}\int_{\T^n} u^\ep_t(x,t)\sig^{\ep,z}(x,t)\,dx=0.
    \]
    Therefore,
    \begin{align*}
        u^\ep_t(z,T) &=\frac{1}{T}\int_0^T\int_{\T^n}u^\ep_t\sig^{\ep,z}\,dxdt\\
        &=\frac{1}{T}\int_0^T\int_{\T^n}\left(-H(x,Du^\ep)+\ep \Delta u^\ep\right)\sig^{\ep,z}\,dxdt.
    \end{align*}
    
\end{proof}

The following lemma is similar to Lemma \ref{lem:NAM}.
\begin{lem}\label{lem:LT2}
    Assume {\rm (A1)--(A2)}.
    There exists $C>0$ independent of $\ep\in (0,1)$ and $T>0$ such that
    \[
    \ep \int_0^T \int_{\T^n} |D^2u^\ep|^2 \sig^{\ep,z}\,dxdt \leq C(1+T).
    \]
    In particular, for $T\geq 1$,
    \[
    \frac{1}{T}\int_0^T \int_{\T^n}\ep |\Delta u^\ep| \sig^{\ep,z}\,dxdt \leq C \ep^{1/2}.
    \]
\end{lem}

\begin{proof}
    The proof of the first inequality was already given in the proof of Lemma \ref{lem:NAM}.

    The second inequality follows from the usual Holder inequality
    \begin{align*}
        &\frac{1}{T}\int_0^T \int_{\T^n}\ep |\Delta u^\ep| \sig^{\ep,z}\,dxdt \\
        \leq \ & \frac{1}{T}\left(\int_0^T \int_{\T^n}\ep^2 |\Delta u^\ep|^2 \sig^{\ep,z}\,dxdt\right)^{1/2}\left(\int_0^T \int_{\T^n} \sig^{\ep,z}\,dxdt\right)^{1/2}\\
        \leq \ & \frac{C}{T}\left(\int_0^T \int_{\T^n}\ep^2 |D^2 u^\ep|^2 \sig^{\ep,z}\,dxdt\right)^{1/2}\left(\int_0^T \int_{\T^n} \sig^{\ep,z}\,dxdt\right)^{1/2}\\
        \leq \ &\frac{C((1+T)\ep)^{1/2}\,T^{1/2}}{T}
        \leq C \ep^{1/2}.
    \end{align*}
\end{proof}

\begin{lem}\label{lem:LT3}
    Assume {\rm (A1)--(A2)} and $\ol H(0)=0$.
    Then, there exists $C>0$ independent of $\ep\in (0,1)$ and $T>0$ such that
    \[
    u^\ep_t(z,T) \geq -\frac{C}{T}- C\ep^{1/2}.
    \]
\end{lem}

\begin{proof}
    Thanks to Lemmas \ref{lem:LT1}--\ref{lem:LT2},
    \begin{align*}
     u^\ep_t(z,T) &= \frac{1}{T}\int_0^T\int_{\T^n}\left(-H(x,Du^\ep)+\ep \Delta u^\ep\right)\sig^{\ep,z}\,dxdt\\
     &\geq \frac{1}{T}\int_0^T\int_{\T^n}-H(x,Du^\ep)\sig^{\ep,z}\,dxdt - C\ep^{1/2}.
    \end{align*}

    On the other hand, the formula of $u^\ep(z,T)$ in Theorem \ref{thm:rep1} gives
    \begin{align*}
    &u^\ep(z,T) \\
    = &\int_{\T^n} g(x)\sigma^{\ep,z}(x,0)\,dx + \int_0^T\int_{\T^n} \left(D_pH(x,Du^\ep)\cdot Du^\ep - H(x,Du^\ep) \right)\sig^{\ep,z}\,dxdt.
    \end{align*}
    As $v_0-\|v_0\|_{L^\infty(\T^n)}-\|g\|_{L^\infty(\T^n)} \leq u \leq v_0+\|v_0\|_{L^\infty(\T^n)}+\|g\|_{L^\infty(\T^n)}$,
    \[
    \left|\int_0^T\int_{\T^n} \left(D_pH(x,Du^\ep)\cdot Du^\ep - H(x,Du^\ep) \right)\sig^{\ep,z}\,dxdt\right| \leq 2 \left(\|v_0\|_{L^\infty(\T^n)}+\|g\|_{L^\infty(\T^n)}\right).
    \]
    Thanks to (A2),
    \[
    D_pH(x,Du^\ep)\cdot Du^\ep - H(x,Du^\ep) \geq \theta H(x,Du^\ep),
    \]
    which implies
    \begin{equation}\label{eq:LT2}
        \int_0^T\int_{\T^n}  H(x,Du^\ep)\sig^{\ep,z}\,dxdt \leq \frac{2 \left(\|v_0\|_{L^\infty(\T^n)}+\|g\|_{L^\infty(\T^n)}\right)}{\theta}.
    \end{equation}
    Thus,
    \begin{equation}\label{eq:LT3}
        u^\ep_t(z,T) \geq -\frac{1}{T}\int_0^T\int_{\T^n} H(x,Du^\ep)\sig^{\ep,z}\,dxdt - C\ep^{1/2} \geq -\frac{C}{T}- C\ep^{1/2}.
    \end{equation}
\end{proof}

\begin{cor}\label{cor:LT1}
    Assume {\rm (A1)--(A2)} and $\ol H(0)=0$.
    Let $\mu^{z,T}$ be the measure as in the statement of Theorem \ref{thm:rep1}.
    Then, for a.e. $(z,T)\in \T^n \times (0,\infty)$,
    \[
   u_t(z,T)=- \frac{1}{T} \int_{\T^n\times \R^n\times [0,T]} H(x,p)\,d\mu^{z,T}(x,p,t).
    \]
    Further, there exists $C>0$ independent of $T>0$ such that
    \begin{equation}\label{eq:u-t}
    u_t(z,T) \geq -\frac{C}{T} \quad \text{ in } \T^n \times (0,\infty)
    \end{equation}
    in the viscosity sense.
\end{cor}

\begin{proof}
    The first equality follows directly from Lemmas \ref{lem:LT1}--\ref{lem:LT2}.

    \smallskip
    
    Secondly, by the stability of viscosity solutions, the fact that $u^\ep \to u$ locally uniformly on $\T^n\times [0,\infty)$, and \eqref{eq:LT3}, we get the inequality \eqref{eq:u-t} in the viscosity sense.
    Of course, as \eqref{eq:u-t} is linear, it also holds in the almost everywhere sense.
\end{proof}

We are now ready to prove the large time behavior result.
\begin{proof}[Proof of Theorem \ref{thm:large time}]
    For each $\ep\in (0,1)$, let $T_\ep=\ep^{-1/4}$.
    By Lemma \ref{lem:LT3}, we have
    \[
    u^\ep_t(x,T_\ep)\geq -C \ep^{1/4} - C\ep^{1/2} \geq - C\ep^{1/4}.
    \]
    Hence,
    \begin{equation}\label{eq:LT4}
    H(x,Du^\ep(x,T_\ep)) - \ep \Delta u^\ep =-u^\ep_t(x,T_\ep)\leq C \ep^{1/4} \quad \text{ in } \T^n.
    \end{equation}
    Further, by Remark \ref{rem:1} (or Lemma \ref{lem:rate 1-2} below),
    \begin{equation}\label{eq:LT5}
    |u^\ep(x,T_\ep) - u(x,T_\ep)| \leq C(1+T_\ep) \ep^{1/2} \leq C \ep^{1/4}.
    \end{equation}
    As $\|Du\|_{L^\infty(\T^n\times [0,\infty))} \leq C$, by the Arzel\`a--Ascoli theorem, we can find a sequence $\{\ep_k\} \to 0$ such that
    \[
    \lim_{k\to\infty}\|u(\cdot,T_{\ep_k}) - v\|_{L^\infty(\T^n)}=0
    \]
    for some $v\in \Lip(\T^n)$.
    Combining this with \eqref{eq:LT5}, we see that
     \[
    \lim_{k\to\infty}\|u^{\ep_k}(\cdot,T_{\ep_k}) - v\|_{L^\infty(\T^n)}=0
    \]
    Write $T_k=T_{\ep_k}$ for $k\in \N$ for simplicity.
    By the stability result for viscosity subsolutions of \eqref{eq:LT4}, we get that $v$ solves
    \begin{equation}\label{eq:LT6}
        H(x,Dv) \leq 0 \quad \text{ in } \T^n.
    \end{equation}
    In particular, $v-C$ is a stationary subsolution to \eqref{eq:HJ} for any constant $C\in \R$.

    \medskip
    
    We now show that $u(\cdot,t)$ converges to $v$ uniformly on $\T^n$ as $t\to \infty$.
    Assume by contradiction that this is not the case.
    Then, we can find a sequence $\{t_k\} \to \infty$ and $w\in \Lip(\T^n)$, $w \neq v$ such that
    \[
    \lim_{k \to \infty} \|u(\cdot,t_k) - w\|_{L^\infty(\T^n)}=0.
    \]
    We note that
    \[
    v(x) - \|u(\cdot,T_k)-v\|_{L^\infty(\T^n)} \leq u(\cdot,T_k).
    \]
    As $v(x) - \|u(\cdot,T_k)-v\|_{L^\infty(\T^n)}$ is a stationary subsolution to \eqref{eq:HJ}, and $u(\cdot,T_k+\cdot)$ is the solution to \eqref{eq:HJ} with initial data $u(\cdot,T_k)$, the comparison principle implies that, for $t \geq T_k$,
    \[
    v - \|u(\cdot,T_k)-v\|_{L^\infty(\T^n)} \leq u(\cdot,T_k+(t-T_k))=u(\cdot,t).
    \]
    Since $\lim_{k\to\infty}\|u(\cdot,T_k)-v\|_{L^\infty(\T^n)}=0$, we deduce
    \[
    v\leq \lim_{k\to \infty}u(\cdot, t_k)=w.
    \]
    By a similar logic, $w\leq v$, and thus, $w=v$, which is a contradiction.
    Therefore,
    \[
    \lim_{t\to\infty}\|u(\cdot,t)-v\|_{L^\infty(\T^n)}=0.
    \]
\end{proof}

\begin{rem}\label{rem:weaker A2}
Some comments are in order.
\begin{itemize}
\item[(i)] As $\|Du^\ep\|_{L^\infty(\T^n\times [0,\infty))} \leq C$, we actually only need to require (A2) holds for $(x,p)\in \T^n\times B(0,C)$ in Theorem \ref{thm:large time}, that is, for $(x,p)\in \T^n\times B(0,C)$,
\[
    D_pH(x,p)\cdot p  \geq (\theta+1) H(x,p) .
\]
The behavior of $H(x,p)$ for $x\in \T^n$ and $|p| \geq C$ is irrelevant.

\item[(ii)] It is clear that (A2) can be replaced by the following assumption, which gives the upper bounds for $u^\ep_t$ and $u_t$.

\smallskip

\noindent (A3) There exists $\theta>0$ such that, 
for $(x,p)\in \T^n\times \R^n$,
\[
    D_pH(x,p)\cdot p  \leq \left(\theta+1\right) H(x,p) .
\]
    
\item[(iii)] We believe that the result in Corollary \ref{cor:LT1} is new and of independent interest.
Moreover, Theorem \ref{thm:large time} can also be proved directly using Corollary \ref{cor:LT1}.
\end{itemize}
\end{rem}

\subsection{Further properties of the solution}
We have a bit more understanding on the behavior of $u$.

\begin{lem}\label{lem:LT-more}
    Assume {\rm (A1)--(A2)}.
    Let $\sig^{z,T}$ and $\mu^{z,T}$ be the measures as in the statement of Theorem \ref{thm:rep1}.
    Then, for a.e. $(z,T)\in \T^n \times (0,\infty)$,
    \[
    u(z,T)+\theta T u_t(z,T) \geq \int_{\T^n} g(x)\,d\sig^{z,T}.
    \]
    In particular, if $g\geq 0$, then $u \geq 0$, and for $z\in \T^n$,
    \begin{equation}\label{eq:u increase}
    T \mapsto T^{1/\theta} u(z,T)  \text{ is nondecreasing in } (0,\infty).
    \end{equation}
   
\end{lem}

\begin{proof}
    By Theorem \ref{thm:rep1} and Corollary \ref{cor:LT1}, we have, for a.e. $(z,T)\in \T^n \times (0,\infty)$,
    \begin{align*}
         u(z,T)
     &= \int_{\T^n\times \R^n\times [0,T]} \left(D_pH(x,p)\cdot p - H(x,p) \right) \,d\mu^{z,T}(x,p,t)+\int_{\T^n} g(x)\,d\sigma^{z,T}(x)\\
     &\geq \theta \int_{\T^n\times \R^n\times [0,T]} H(x,p)  \,d\mu^{z,T}(x,p,t)+\int_{\T^n} g(x)\,d\sigma^{z,T}(x)\\
     &=-\theta T u_t(z,T)+\int_{\T^n} g(x)\,d\sigma^{z,T}(x).
    \end{align*}
     Hence,
    \[
    u(z,T)+\theta T u_t(z,T) \geq \int_{\T^n} g(x)\,d\sig^{z,T}.
    \]
    
    If $g \geq 0$, then we have
    \[
    u(z,T)+\theta T u_t(z,T) \geq 0 \qquad \text{ in } \T^n \times (0,\infty)
    \]
    in the viscosity sense.
    In particular,
    \[
    \frac{d}{dT}\left( T^{1/\theta}u(z,T) \right) = \frac{T^{1/\theta-1}}{\theta}(u(z,T)+\theta T u_t(z,T)) \geq 0.
    \]
    We conclude that $T^{1/\theta}u(z,T) \geq 0$ and $ T \mapsto T^{1/\theta} u(z,T) $   is nondecreasing in  $(0,\infty)$.
    
\end{proof}

\begin{rem}
    Let $p=0$ in (A2) to get
    \[
    H(x,0) \leq 0 \quad \text{ for } x\in \T^n.
    \]
    In particular, the constant function $\min_{\T^n}g$ is a subsolution to \eqref{eq:HJ}, and thus, by the usual comparison principle,
    \[
    u(x,t) \geq \min_{\T^n} g \qquad \text{ for } (x,t)\in \T^n\times [0,\infty).
    \]
    This gives another proof of a claim in Lemma \ref{lem:LT-more} above that $u \geq 0$ if $g\geq 0$.
    It is not clear to us yet how useful is the nondecreasing property of $ T \mapsto T^{1/\theta} u(z,T) $.
\end{rem}

For completeness, we now show the $O(\ep^{1/2})$ convergence rate of $u^\ep$ to $u$ via the nonlinear adjoint method, which was already done in \cite{Evans1, Tran1, Cagnetti-Gomes-Tran, CGMT, Le-Mitake-Tran, tranbook, CGM23, Cirant-Goffi}.
See \cite{Guo-Tran-Zhang} for policy iteration and related problems for nonconvex viscous Hamilton--Jacobi equations.

\begin{lem}\label{lem:rate 1-2}
Assume {\rm (A1)}.
Then, there exists $C>0$ independent of $\ep\in (0,1)$ such that, for $(z,T)\in \T^n\times [0,\infty)$ and $\ep\in (0,1)$,
\[
|u^\ep(z,T)-u(z,T)| \leq C(1+T) \ep^{1/2}.
\]
\end{lem}

\begin{proof}
    We note that $u^\ep$ is smooth enough with respect to $\ep\in (0,1)$.
    Denote by $u^\ep_\ep = \frac{\partial u^\ep}{\partial \ep}$.
    As $u^\ep(x,0)=g(x)$ for all $\ep\in (0,1)$, $u^\ep_{\ep}(x,0)=0$ for $x\in \T^n$.

    \smallskip
    
    Differentiate \eqref{eq:HJ-ep} with respect to $\ep$ to get
    \[
    \cL^\ep[u^\ep_\ep]=(u^\ep_\ep)_t +D_pH(x,Du^\ep)\cdot Du^\ep_\ep - \ep \Delta u^\ep_\ep=\Delta u^\ep.
    \]
    
    Multiply this equality by $\sig^{\ep,z}$ and integrate to imply
    \begin{align*}
    |u^\ep_\ep(z,T)|&=\left|\int_0^T\int_{\T^n}\Delta u^\ep \sig^{\ep,z}\,dxdt\right|\\
    &\leq C\left(\int_0^T\int_{\T^n}|D^2 u^\ep|^2 \sig^{\ep,z}\,dxdt\right)^{1/2}\left(\int_0^T\int_{\T^n} \sig^{\ep,z}\,dxdt\right)^{1/2}\\
    &\leq C\ep^{-1/2}(1+T)^{1/2} \, T^{1/2} \leq C(1+T) \ep^{-1/2}.
    \end{align*}
   We used Lemma \ref{lem:LT2} in the second last inequality above.
   Then, by the fundamental theorem of calculus,
   \[
   |u^\ep(z,T)-u(z,T)| \leq \int_0^\ep \left|\frac{\partial u^\delta}{\partial\delta}(z,T)\right|\,d\delta
   \leq C(1+T)\int_0^\ep \delta^{-1/2}\,d\delta \leq C(1+T) \ep^{1/2}.
   \]
   As discussed in \cite{QSTY}, this $O(\ep^{1/2})$ convergence rate is optimal for general Hamiltonians satisfying (A1).
\end{proof}

\subsection{Examples and discussions}\label{subsec:LT-A2}
Let us first discuss assumption (A2).
Fix $p\in \R^n$, and define, for $s\geq 0$,
\[
\varphi_p(s)=H(x,sp).
\]
Then, (A2) holds if and only if for $p\in \R^n$ and $s>0$,
\[
\varphi_p'(s) = D_pH(x,sp)\cdot p=\frac{1}{s} D_pH(x,sp)\cdot (sp) \geq \frac{1}{s}(\theta+1)H(x,sp) = \frac{\theta+1}{s}\varphi_p(s).
\]
This happens if and only if
\[
\left(s^{-(\theta+1)} \varphi_p(s)\right)' \geq 0.
\]
Thus, (A2) is equivalent to the following assumption.
\begin{itemize}
    \item[(A2’)] There exists $\theta>0$ such that:  For $(x,p)\in \T^n\times \R^n$,
    \[
    s \mapsto s^{-(\theta+1)}H(x,sp) \text{ is nondecreasing in } (0,\infty).
    \]
\end{itemize}
In particular, this shows that assumption (A2) (or (A2’)) is about the behavior of $H(x,p)$ and its directional derivative along the $p$ direction, and not about the full gradient or strict convexity.
This is demonstrated clearly in Examples \ref{ex:H2} and \ref{ex5} below, where conditions (A6)$_\pm$ in \cite{BIM} and (H4) in \cite{BS1} do not hold in general.
As noted, Theorem \ref{thm:large time} complements the results in \cite{BS1, BIM}, which pushes further the study of large time behavior of nonconvex first-order Hamilton--Jacobi equations.

\medskip

Before getting to the examples, we give a different proof of \eqref{eq:u increase} by using (A2').

\begin{proof}[Another proof of \eqref{eq:u increase}]
    For $\lambda \geq 1$, define
    \[
    v(x,t) = \lambda^{1/\theta}u(x,\lambda t) \qquad \text{ for } (x,t)\in \T^n\times [0,\infty).
    \]
    As $g\geq 0$ and $\lambda \geq 1$, for $x\in \T^n$,
    \[
    v(x,0)=\lambda^{1/\theta}u(x,0)=\lambda^{1/\theta} g(x) \geq g(x).
    \]
    Thanks to (A2'), 
    \begin{align*}
        v_t(x,t) + H(x,Dv(x,t))
        &=  \lambda^{(\theta+1)/\theta} u_t(x,\lambda t) + H(x, \lambda^{1/\theta}Du(x,\lambda t))\\
        &\geq \lambda^{(\theta+1)/\theta} u_t(x,\lambda t) +\lambda^{(\theta+1)/\theta} H(x, Du(x,\lambda t))=0.
    \end{align*}
    Hence, $v$ is a supersolution to \eqref{eq:HJ}.
    By the usual comparison principle, for $(x,t)\in \T^n\times [0,\infty)$,
    \[
    u(x,t) \leq v(x,t) = \lambda^{1/\theta}u(x,\lambda t),
    \]
    which implies \eqref{eq:u increase}.
\end{proof}

Let us now give some examples of $H$ satisfying (A1)--(A2) and $\ol H(0)=0$.
\begin{ex}\label{ex:H1}
    Let
    \[
    H(x,p)=H(p)=|p|^4-|p|^2 \qquad \text{ for } (x,p) \in \T^n\times \R^n.
    \]
    Then, $\ol H(0)=0$ as $H(0)=0$.
    It is clear that {\rm (A1)} holds.
    We compute
    \begin{align*}
        &D_pH(x,p)\cdot p - H(x,p)=D_pH(p)\cdot p - H(p)\\
        =\ &(4|p|^2p - 2p)\cdot p - (|p|^4-|p|^2)=3|p|^4-|p|^2 \geq H(p).
    \end{align*}
    Hence, {\rm (A2)} holds with $\theta=1$.
    Besides,
    \[
    D^2 H(p)=(4|p|^2-2)I_n + 8p\otimes p .
    \]
    Of course, $H$ is neither concave or convex. 
    See Figure \ref{fig:H1}.
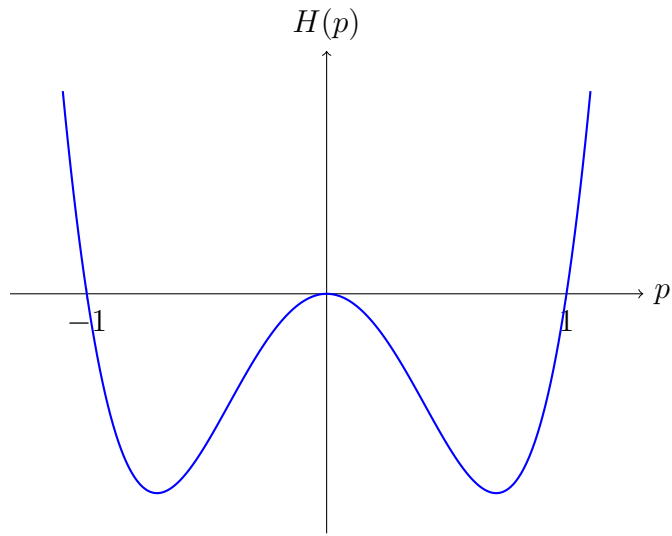
\begin{figure}[htp]
\begin{center}
\begin{tikzpicture}
  \begin{axis}[
      axis lines = middle,
      xlabel = $p$,
      ylabel = {$H(p)$},
      domain = -1.1:1.1,
      samples = 200,
      width=10cm,
      height=8cm,
      enlargelimits=true,
      xtick={-2,-1,0,1,2},
      ytick={-1,-0.5,0,0.5,1},
      axis line style={->},
      xlabel style={right},
      ylabel style={above},
    ]
    \addplot[blue, thick] {x^4 - x^2};
  \end{axis}
\end{tikzpicture}
\caption{Graph of $H(p)=|p|^4-|p|^2$ in one dimension}\label{fig:H1}
\end{center}
\end{figure}
\end{ex}

\begin{ex}\label{ex:H2}
    Let
    \[
    H(x,p)=c(x)(a(p)- b(p)) \qquad \text{ for } (x,p) \in \T^n\times \R^n.
    \]
    Here, $c\in C^2(\T^n,(0,\infty))$, $a,b\in C^2(\R^n, [0,\infty)) \cap C^2(\R^n\setminus\{0\}, (0,\infty))$ are given such that $a$ is homogeneous of degree $4$, and $b$ is homogeneous of degree $2$.
    As $a(0)=b(0)=0$, $H(x,0)=0$ for $x\in \T^n$, and so $\ol H(0)=0$.
    
    It is clear that {\rm (A1)} holds as for $p\neq 0$,
    \[
    H(x,p)=c(x)\left(a\left(\frac{p}{|p|}\right)|p|^4-b\left(\frac{p}{|p|}\right)|p|^2\right).
    \]
    Set $\theta=1$.
    For $p \in \R^n$ and $s>0$, 
    \[
    s^{-2}H(sp)=s^2 c(x)a(p) - c(x)b(p),
    \]
    which is nondecreasing for $s\in (0,\infty)$.
    Therefore, {\rm (A2’)} and {\rm (A2)} hold with $\theta=1$.
    In particular, for each $x\in \T^n$, the $0$-sublevel set of $H(x,\cdot)$ is
    \begin{align*}
    S_x&=\{p\,:\,H(x,p) \leq 0\}\\
     &=\{0\} \bigcup \left\{p \neq 0\,:\,|p| \leq a\left(\frac{p}{|p|}\right)^{-1/2}b\left(\frac{p}{|p|}\right)^{1/2}  \right\}.
    \end{align*}
    Generally speaking, $S_x$ is not convex, and thus, conditions {\rm (A6)$_\pm$} in \cite{BIM} do not hold.
    We note that $H$ given here is just a prototypical one, and there are many similar Hamiltonians satisfying {\rm(A1)--(A2)} and $\ol H(0)=0$.
\end{ex}

\medskip

Let us now recall an example in \cite{BS1} where large time behavior fails if $H$ is merely convex in $p$.

\begin{ex}\label{ex:H3}
    Assume $n=1$ and
    \[
    H(x,p)=H(p)=|p+10|-10.
    \]
    Consider
    \[
    \begin{cases}
        u_t+|u_x+10|-10 = 0 \qquad &\text{ in } \T\times (0,\infty),\\
        u(x,0)=\sin(2\pi x) \qquad &\text{ on } \T.
    \end{cases}
    \]
    Then, $u(x,t)=\sin(2\pi(x-t))$ for $(x,t)\in \T\times [0,\infty)$ is the solution to the above, and $u(x,t)$ does not converge as $t\to\infty$.
\end{ex}

Basically, Example \ref{ex:H3} shows that we could have traveling wave solutions, which do not converge to a solution of the corresponding cell problem, if $H(x,p)=H(p)$ is linear in a neighborhood of $p=0$.

\medskip

As noted, \cite{BS1, BIM} obtained large time behavior for \eqref{eq:HJ} for some possibly nonconvex Hamiltonians.
    Roughly speaking, when $\ol H(0)=0$, \cite{BS1, BIM} only need to require some conditions similar to the strict convexity of $H$ near the $0$-sublevel set of $H$ (conditions (A6)$_\pm$ in \cite{BIM} or condition (H4) in \cite{BS1}).
    These are rather different from assumption (A2) that we put here.
    Let us analyze (H4) in \cite{BS1} to see it a bit more.
    For simplicity, we consider the case
    \begin{equation}\label{eq:BS}
    \begin{cases}
        H(x,p) = K(p) - V(x)  \text{ for } (x,p)\in \T^n\times \R^n,\\
        \text{$K\in C(\R^n)$ is coercive and $K(0)=0=\min_{\R^n} K$},\\
        \text{$V \in C(\T^n)$, $\min_{\T^n} V = 0$, and $\max_{\T^n} V = M$ for some $M>0$.}
    \end{cases}
    \end{equation}
\begin{lem} \label{lem:c-0}
Assume \eqref{eq:BS}. 
Then $\ol H(0)=0$.
\end{lem}

\begin{proof}
Let $w$ be a solution to \eqref{eq:cell}. 
By a priori estimates, $w \in \Lip(\T^n)$ and $w$ solves \eqref{eq:cell} in the almost everywhere sense.

Pick $x_1 \in \T^n$ such that $w(x_1) = \min_{\T^n} w$.
By the viscosity supersolution test,
\[
\ol H(0) \leq H(x_1,0) = - V(x_1) \leq 0.
\]
We need to show the reverse inequality. 
Pick $x_2 \in \T^n$ such that $V(x_2) = \min V = 0$.
Since $w$ solves \eqref{eq:cell} in the almost everywhere sense, there exists a sequence $\{y_k\} \subset \T^n$
such that $\lim_{k \to \infty} y_k =x_2$, $w$ is differentiable at $y_k$, and
\[
\ol H(0)=K(Dw(y_k)) - V(y_k) \geq -V(y_k)
\]
in the classical sense.
Let $k \to \infty$ in the above to get the desired conclusion.

\end{proof}

Let us now recall condition (H4) in \cite{BS1} with $\ol H(0)=0$ as verified in Lemma \ref{lem:c-0}.
\begin{itemize}
\item[(H4)] There exists $\eta_0>0$ such that, for each $\eta \in (0,\eta_0)$, there exists a constant $\psi_\eta>0$ satisfying:
If $H(x,p) \leq 0$ and $H(x,p+q) \geq \eta$ for some $x\in \T^n$, $p,q \in \R^n$, then
\[
H(x,p+sq) \geq s H(x,p+q) + \psi_\eta (s-1) \qquad \text{ for all } s \geq 1.
\]
\end{itemize}
Denote by
\[
A= \{ p \in \R^n\,:\,  0 \leq K(p) \leq M\}.
\]
We have the following result.
\begin{prop}\label{prop:BS}
Assume \eqref{eq:BS} and {\rm (H4)}.
Then, the following properties hold.
\begin{itemize}
    \item[(i)] $K$ is convex on $A$.
    \item[(ii)] For $p\in A$ and $r\in \R^n \setminus A$,
    \[
    K\left(\frac{p+r}{2} \right) < \frac{K(p)+K(r)}{2}.
    \]
\end{itemize}

\end{prop}

Before providing the proof of this proposition, we need the following result.

\begin{lem} \label{lem:A-conv}
Assume \eqref{eq:BS} and {\rm (H4)}.
Then $A$ is convex.
\end{lem}

\begin{proof}
Assume by contradiction that $A$ is not convex.
Then there exist $p,q \in \R^n$ and $s>1$ such that
\[
p,\, p +sq \in A, \qquad \text{ and } \qquad p+q \notin A.
\]
Pick $x_0 \in \T^n$ such that $V(x_0) = M = \max_{\T^n} V$.
Then 
\[
H(x_0,p) = K(p) - V(x_0) \leq 0 \quad \text{ and } \quad
H(x_0,p+sq) = K(p+sq) - V(x_0) \leq 0,
\]
but 
\[
H(x_0,p+q) = K(p+q) - M =: \eta>0.
\]
We hence can use (H4) to deduce that
\[
H(x_0, p+sq) \geq s H(x_0, p+q) + \psi_\eta ( s-1) = s \eta + \psi_\eta (s-1) >0,
\]
which is absurd.
Therefore, $A$ is convex.
\end{proof}

\begin{proof}[Proof of Proposition \ref{prop:BS}]
We first prove (i).
Fix $p,r \in A$. 
By Lemma \ref{lem:A-conv}, $(p+r)/2 \in A$ as well. 
We need to show that
\begin{equation}\label{eq:conv}
K\left(\frac{p+r}{2} \right) \leq \frac{K(p)+K(r)}{2}.
\end{equation}
Of course, if 
\[
K\left(\frac{p+r}{2} \right) \leq \min\{K(p),K(r)\},
\]
then \eqref{eq:conv} holds immediately.
Without loss of generality, we only need to consider the case that
\begin{equation}\label{eq:strict convex-1}
K(p) < K\left(\frac{p+r}{2} \right).
\end{equation}
Let $q=(r-p)/2$. 
Pick $y\in \T^n$ such that $V(y)=K(p)$. 
Then
\[
H(y,p)= K(p)- V(y)= 0 \quad \text{and} \quad H(y,p+q) = K(p+q) - V(y)=:\eta>0.
\]
We use condition (H4) with $s=2$ to get
\[
H(y,r) = H(y,p+2q) \geq 2 H(y,p+q) + \psi_\eta = 2 K(p+q) - 2 V(y) + \psi_\eta,
\]
which, together with the fact that $V(y)= K(p)$, yields
\[
K(r) + K(p) = K(p+2q) + K(p) \geq 2 K(p+q) + \psi_\eta = 2 K\left(\frac{p+r}{2} \right) +\psi_\eta.
\]
We conclude that $K$ is convex in $A$.

Next, we prove (ii).
Fix $p\in A$ and $r\in \R^n \setminus A$.
Then, 
\[
K(p) \leq M < K(r).
\]
If 
\[
K\left(\frac{p+r}{2} \right) \leq K(p),
\]
then we are done.
Otherwise, \eqref{eq:strict convex-1} holds and we just repeat the corresponding steps in the proof of (i) to get the desired conclusion.
\end{proof}

\begin{rem}
We have a few further comments.

(i) Based on the above proof, we actually have that, for a line segment $l \subset A$, if $\min_{p \in l} K(p)$ is achieved at exactly one point $p_l \in l$,
then $K$ is strictly convex on $l$.

(ii) Of course, (i) does not imply that $K$ itself is strictly convex on $A$, and in particular,
it does not rule out a possibility that level sets of $K$ contain some line segments.

(iii) The above proof also implies that $K$ is strictly convex near $\partial A$ along the directions emanating from A to $\R^n\setminus A$.
\end{rem}

We now give an example of $H$ satisfying \eqref{eq:BS} and (A2) but not (H4).

\begin{ex}\label{ex5}
    Assume
    \[
     H(x,p) = K(p) - V(x) \qquad  \text{ for } (x,p)\in \T^n\times \R^n,
     \]
    where $K\in C^2(\R^n,[0,\infty))$ is coercive and positive homogeneous of degree $m$ for $m\geq 2$ given.
    Here, $V \in C^2(\T^n)$, $\min_{\T^n} V = 0$, and $\max_{\T^n} V = M$ for some $M>0$.
    Choose $K$ such that $A=\{p\in \R^n\,:\, 0 \leq K(p)\leq M\}$ is not convex.
    We get that \eqref{eq:BS} holds and {\rm (H4)} does not hold.

    For fixed $(x,p)\in \T^\times \R^n$ and $s>0$,
    \[
    s^{-m}H(x,sp) = s^{-m}(K(sp)- V(x))= K(p) - s^{-m} V(x),
    \]
    which is nondecreasing in $s$.
    Hence, we have {\rm (A2')}, which implies that {\rm (A2)} holds.
    We note that this Hamiltonian does not satisfy the conditions in \cite{N-R} too.
\end{ex}
    
\section{Open problems}\label{sec:open problems}

In this section, we list some open problems/questions that might be of interests to the reader.
Firstly, we discuss the measure $\mu^{\ep,z}$, which encodes the information of the characteristics going through $(z,T)$, from the representation formulas in Theorem \ref{thm:rep1}.
\begin{quest}\label{quest1}
    Assume the settings of Theorem \ref{thm:rep1}.
    Show that there are situations where $\{\mu^{\ep,z}\}$ has two different subsequential limits weakly in the sense of measures as $\ep \to 0$.
\end{quest}

\begin{quest}\label{quest2}
    Assume the settings of Theorem \ref{thm:rep1}.
    Assume further that $p\mapsto H(x,p)$ is convex for each $x\in \T^n$.
    Does $\{\mu^{\ep,z}\}$ converge weakly in the sense of measures as $\ep \to 0$ to a unique limit?
\end{quest}
Recall that in the convex setting, we have the optimal control formula \eqref{eq:OC}.
In many situations, \eqref{eq:OC} admits more than one minimizing curves.
For each minimizing curve $\xi$, we showed that it corresponds to a minimizing measure $\gam^\xi$ of \eqref{eq:formula-minimizing-u}.
Further, convex combinations of minimizers of \eqref{eq:formula-minimizing-u} is also a minimizer.
Hence, in general, \eqref{eq:formula-minimizing-u} has infinitely many solutions, and this makes Question \ref{quest2} quite delicate.
Of course, it is also natural to study Questions \ref{quest1}--\ref{quest2} when $T$ is small, that is, before the characteristics have crossed.
\medskip

Next, we focus on the vanishing viscosity approximation of the cell problem \eqref{eq:cell-ep} and the Mather measures.

\begin{quest}\label{quest3}
    Assume {\rm (A1)}.
    For $\ep\in (0,1)$, let $v^\ep$ be the solution to \eqref{eq:cell-ep} with $v^\ep(0)=0$.
    Do we have that $v^\ep$ converges uniformly in $\T^n$ to a unique limit as $\ep \to 0$?
\end{quest}
It seems that the answer for Question \ref{quest3} should be negative, but we do not know very clearly at this moment.
Let us now impose further the uniformly convexity of the Hamiltonian in this question to make it more specific.

\begin{quest}\label{quest4}
    Assume {\rm (A1)} and $H$ is uniformly convex in $p$.
    For $\ep\in (0,1)$, let $v^\ep$ be the solution to \eqref{eq:cell-ep} with $v^\ep(0)=0$.
    Do we have that $v^\ep$ converges uniformly in $\T^n$ to a unique limit as $\ep \to 0$?
\end{quest}

Question \ref{quest4} was already posed in \cite[Section 6.6.2]{Le-Mitake-Tran}, but it is worth to repeat it here.
This question was confirmed positively for some very special cases of $H$ in \cite{AIPS, Bessi}.
Next is a question about Mather measures.

\begin{quest}\label{quest5}
    What can we say about Mather measures if we assume both {\rm (A1)--(A2)} and $\ol H(0)=0$?
    Do we have that the Mather measures are invariant under the Hamiltonian flow?
\end{quest}

It is also natural to ask whether the Mather measures play the role of the uniqueness set for the cell problem \eqref{eq:cell}.

\begin{quest}\label{quest12}
Assume both {\rm (A1)--(A2)} and $\ol H(0)=0$.
Do we have that the Mather measures play the role as the uniqueness set for \eqref{eq:cell}?
\end{quest}
See \cite{Mitake-Tran} and the references therein for this uniqueness property in the convex setting.
Next, we pose some questions concerning the large time behavior and large time limits for the solution $u$ of \eqref{eq:HJ}.

\begin{quest}\label{quest6}
Assume {\rm (A1)}.
Assume further that for each $x\in \T^n$, $p\mapsto H(x,p)$ is not affine on any line segment on $\R^n$.
Do we have the large time behavior result for the solution $u$ of \eqref{eq:HJ}?
\end{quest}

Question \ref{quest6} is quite open-ended. The additional assumption that $p\mapsto H(x,p)$ is not affine on any line segment on $\R^n$ is simply to rule out traveling wave solutions given in Examples \ref{ex:H3}.

The following two questions are related to the next steps of Theorem \ref{thm:large time}.
\begin{quest}\label{quest7}
Assume {\rm (A1)--(A2)} and $\ol H(0)=0$.
Then, we have the large time behavior result in Theorem \ref{thm:large time}.
Characterize the large time limit $v(x)=\lim_{t\to\infty} u(x,t)$ in terms of the initial data $g$ and the Hamiltonian $H$.
\end{quest}

\begin{quest}\label{quest8}
Assume {\rm (A1)--(A2)} and $\ol H(0)=0$.
Then, we have the large time behavior result in Theorem \ref{thm:large time}.
Is it possible to find a quantitative convergence rate of $u(x,t)$ to $v(x)$ as $t\to\infty$?
\end{quest}
Note that Question \ref{quest8} remains open even in the case where $H$ is uniformly convex in $p$.
It is possible that the convergence rate of $u(x,t)$ to $v(x)$ as $t\to\infty$  can be arbitrarily slow.
See \cite{F-U} for some related results.

\medskip

As noted in the introduction, large time behavior for possibly degenerate viscous Hamilton--Jacobi equations with uniformly convex Hamiltonians was proved in \cite{CGMT}.
There is not yet any result in the literature for nonconvex Hamiltonians, which leads to the following question.

\begin{quest}\label{quest9}
    Study the large time behavior of possibly degenerate viscous Hamilton--Jacobi equations when the Hamiltonian is not convex.
\end{quest}

Along this line, it is important to develop a theory of Mather measures for possibly degenerate viscous Hamilton--Jacobi equations.

\begin{quest}\label{quest10}
    Develop a theory of Mather measures for possibly degenerate viscous Hamilton--Jacobi equations when the Hamiltonian is not convex.

\end{quest}

Finally, we pose a question on large time behavior for fully nonlinear parabolic equations.

\begin{quest}\label{quest11}
    Study the large time behavior of
    \begin{equation*}
    \begin{cases}
        u_t + F(x,Du, D^2 u)=0 \qquad &\text{ in } \T^n \times (0,\infty),\\
        u(x,0)=g(x) \qquad &\text{ on } \T^n.
    \end{cases}
\end{equation*}
Here, $F$ is degenerate elliptic.
\end{quest}

\section*{Acknowledgment}
I would like to thank Professor Hiroyoshi Mitake for his important comments and suggestions.

\end{document}